\documentclass[UTF-8,reqno]{amsart}
\usepackage{enumerate}
\usepackage{tabularx,array}\usepackage{tabulary}
\usepackage{xltabular}\usepackage{xltabular}\usepackage{longtable}\usepackage{lipsum}
\usepackage{verbatim,listings}
\usepackage{amssymb,url,color, booktabs}

\usepackage{mathrsfs}
\usepackage{amsmath}
\usepackage{verbatim}
\usepackage{fancyhdr}
\usepackage[nobysame]{amsrefs}
\BibSpec{article}{%
+{}{\PrintAuthors} {author}
+{,}{ \textrm} {title}
+{.}{ \textit} {journal}
+{,}{ \textbf} {volume}
+{}{ \parenthesize} {date}
+{,}{ } {pages}
+{.}{ arXiv:} {eprint}
+{.}{} {transition}
}
\BibSpec{book}{%
+{}{\PrintAuthors} {author}
+{,}{ \textit} {title}
+{.}{ \textrm} {series} 
+{,}{ Vol.} {volume}
+{.}{ } {publisher}
+{,}{ } {date}
+{.}{} {transition}
}
\usepackage{color}
\usepackage{tikz}
\usepackage{graphicx}
\usepackage[colorlinks=true]{hyperref}
\hypersetup{
    linkcolor=blue,          
    citecolor=red,        
    filecolor=blue,      
    urlcolor=cyan
}

\numberwithin{equation}{section}

\newcommand{\be}{\begin{eqnarray}}
\newcommand{\ee}{\end{eqnarray}}
\newcommand{\ce}{\begin{eqnarray*}}
\newcommand{\de}{\end{eqnarray*}}
\newtheorem{theorem}{Theorem}[section]

\newtheorem{remark}{Remark}[section]

\newtheorem{lemma}[theorem]{Lemma}
\newtheorem{definition}{Definition}[section]
\newtheorem{proposition}[theorem]{Proposition}
\newtheorem{Examples}[theorem]{Example}
\newtheorem{corollary}[theorem]{Corollary}

\newenvironment{proof of theorem 1.4}{{\it Proof of Theorem 1.4}.}{{\hfill 	
		$\square$\hskip - \parfillskip}}
\newenvironment{proof of theorem 1.5}{{\it Proof of Theorem 1.5}.}{{\hfill 	
		$\square$\hskip - \parfillskip}}
\newenvironment{proof of theorem 1.6}{{\it Proof of Theorem 1.6}.}{{\hfill 	
		$\square$\hskip - \parfillskip}}

\makeatletter

\newcommand{\Rmnum}[1]{\expandafter\@slowromancap\romannumeral #1@}
\makeatother

\usepackage{caption}
\def\[{{\Big[}}
\def\]{{\Big]}}
\def\<{{\langle}}
\def\>{{\rangle}}
\def\({{\Big(}}
\def\){{\Big)}}

\def\bx{{\mathbf{x}}}

\def\min{{\mathord{{\rm min}}}}

\def\={&\!\!=\!\!&}

\def\1{{\mathbf{1}}}

\def\geq{\geqslant}
\def\leq{\leqslant}
\def\ge{\geqslant}
\def\le{\leqslant}

\def\[{{\Big[}}
\def\]{{\Big]}}
\def\<{{\langle}}
\def\>{{\rangle}}
\def\({{\Big(}}
\def\){{\Big)}}

\def\bx{{\mathbf{x}}}

\def\W{{\mathcal W}}

\def\min{{\mathord{{\rm min}}}}

\def\={&\!\!=\!\!&}
\def\bt{\begin{theorem}}
\def\et{\end{theorem}}
\def\bl{\begin{lemma}}
\def\el{\end{lemma}}
\def\br{\begin{remark}}
\def\er{\end{remark}}
\def\bx{\begin{Examples}}
\def\ex{\end{Examples}}
\def\bd{\begin{definition}}
\def\ed{\end{definition}}
\def\bp{\begin{proposition}}
\def\ep{\end{proposition}}
\def\bc{\begin{corollary}}
\def\ec{\end{corollary}}

\def\geq{\geqslant}
\def\leq{\leqslant}
\def\ge{\geqslant}
\def\le{\leqslant}

\def\<{\langle} \def\>{\rangle}

\def\bpf{\begin{proof}}
\def\epf{\end{proof}}

\allowdisplaybreaks

	\tikzset{
	pattern size/.store in=\mcSize,
	pattern size = 5pt,
	pattern thickness/.store in=\mcThickness,
	pattern thickness = 0.3pt,
	pattern radius/.store in=\mcRadius,
	pattern radius = 1pt}
\makeatletter
\pgfutil@ifundefined{pgf@pattern@name@_ak5ydppfi}{
	\pgfdeclarepatternformonly[\mcThickness,\mcSize]{_ak5ydppfi}
	{\pgfqpoint{0pt}{0pt}}
	{\pgfpoint{\mcSize+\mcThickness}{\mcSize+\mcThickness}}
	{\pgfpoint{\mcSize}{\mcSize}}
	{
		\pgfsetcolor{\tikz@pattern@color}
		\pgfsetlinewidth{\mcThickness}
		\pgfpathmoveto{\pgfqpoint{0pt}{0pt}}
		\pgfpathlineto{\pgfpoint{\mcSize+\mcThickness}{\mcSize+\mcThickness}}
		\pgfusepath{stroke}
}}
\makeatother
\tikzset{every picture/.style={line width=0.75pt}} 

\begin{document}
\title{Stability of $L^p$-Brunn--Minkowski inequalities for volume and intrinsic volumes}
\author{Shanwei Ding, Jinyu Gao$^{*}$, and Guanghan Li}\thanks{* Corresponding author}


\begin{abstract}
In this paper, we establish stability results of $L^p$-Brunn--Minkowski inequalities as well as $L^p$-Minkowski inequalities for volume and intrinsic volumes in the non-origin-symmetric setting and for the range $p>-k+1-\frac{2k}{n}$. We prove quantitative stability estimates in terms of both the $L^2$ distance and the relative asymmetry index. As a consequence, the corresponding $L^p$-Brunn--Minkowski inequalities,  $L^p$-Minkowski inequalities, and the uniqueness for the $L^p$-Christoffel-Minkowski problem are obtained.
\end{abstract}

\maketitle
\setcounter{tocdepth}{2}
\tableofcontents

\section{Introduction}\label{sec 1}

The Brunn--Minkowski theory is a cornerstone of modern convex geometry and geometric analysis, and it has deep connections with analysis, probability, and partial differential equations. Its development spans more than a century and has had a profound influence across mathematics. At its core lies the classical Brunn--Minkowski inequality (see \cite{book-convex-body,Gardner-book-2006,Gardner-2002,Borell-1997}) and its $L^p$ analogue, the so-called $L^p$-Brunn--Minkowski inequality introduced by \cite{Firey-book} and later developed by  \cite{Lutwak-Yang-Zhang-2012,L93}.

Let $\mathcal{K}_o^{n+1}$ with $n\ge 1$ denote the set of convex bodies (compact, convex sets with non-empty interiors) in $\mathbb{R}^{n+1}$ containing the origin in their interiors.
We write $\mathbb{B}^{n+1}:=\{x\in\mathbb{R}^{n+1}:|x|\leq 1\}$ for the standard unit ball and $\mathbb{S}^{n}:=\partial\mathbb{B}^{n+1}$ for the unit sphere.
For  $K\in\mathcal{K}^{n+1}_o$, its support function $u_K: \mathbb{S}^{n} \rightarrow \mathbb{R}$ is defined by
$$
u_K(x)=\max _{y \in K}\langle x, y\rangle, \quad x \in \mathbb{S}^n,
$$
where $\langle\cdot,\cdot\rangle$ denotes the standard Euclidean inner product on $\mathbb{R}^{n+1}$.

For $K,L\in\mathcal{K}_o^{n+1}$, $p\in\mathbb{R}$, and $t\in(0,1)$, the $L^p$-Minkowski convex combination (or, simply, the $L^p$-Minkowski sum) of $K$ and $L$ is defined by
\begin{equation}
	(1-t) \cdot K+_p t \cdot L:=\left\{\begin{array}{ll}
		W\left(\left((1-t) u_K^p+t u_L^p\right)^{\frac{1}{p}}\right), & p \neq 0, \\
		W\left(u_K^{1-t} u_L^t\right), & p=0.
	\end{array} \right.
\end{equation}
Here the {Wulff shape} $W(f)$ associated with $f$ is the unique maximal element in $$\left\{Q\in\mathcal{K}_o^{n+1}:u_Q\leq f\right\}.$$
For $p=1$, this definition recovers the classical Minkowski convex combination: $(1-t)\cdot K+_1t\cdot L=(1-t)K+tL=\{(1-t)x+ty:x\in K,\ y\in L\}$. In this notation, the $L^p$-Brunn--Minkowski inequality \cite{L93} asserts that, for $p\geq 1$ and $t\in(0,1)$,
\begin{align}\label{equ:BM-neq-p=1}
	V_n\left((1-t) \cdot K+_p t \cdot L\right) \geq\left((1-t) V_n(K)^{\frac{p}{n+1}}+t V_n(L)^{\frac{p}{n+1}}\right)^{\frac{n+1}{p}},
\end{align}
holds with equality if and only if $K$ and $L$ are dilates (and may differ by a translation in the case $p=1$). Here $\frac{1}{n+1}V_n(K)=|K |$ denotes the volume of $K\subset \mathbb{R}^{n+1}$. In the case $p=1$, it reduces to the well-known Brunn--Minkowski inequality.

The substantially more delicate regime $p<1$, which yields stronger forms of the classical Brunn--Minkowski inequality, has attracted considerable attention. A closely related formulation is given by the $L^p$-Minkowski inequality, which is equivalent in the present setting (see \cite{Boroczky-Lutwak-Yang-Zhang-2012,Milman-JEMS-2025}):
\begin{align}
	\label{equ:M-neq}
	\frac{1}{p}\int_{\mathbb{S}^n}u_L^pu_K^{1-p}E_n(\lambda_K)~dx&\geq\frac{1}{p} V_n(K)^{\frac{n+1-p}{n+1}}V_n(L)^{\frac{p}{n+1}},\quad \mathrm{ if }\ p\ne0,\\
	\label{equ:M-neq=p=1}
	\int_{\mathbb{S}^n}\left(\log(\frac{u_L}{u_K})\right)u_KE_n(\lambda_K) ~dx
	&\geq \frac{1}{n+1} V_n(K)\log\frac{V_n(L)}{V_n(K)},\quad \ \! \mathrm{ if }\ p=0.
\end{align}
For developments concerning these inequalities, see, for example, \cite{Boroczky-Lutwak-Yang-Zhang-2012,Chen-Huang-Li-Liu-2020,Putterman-2021,Milman-JEMS-2025,Kolesnikov-Milman-MAMS-2022,Milman-JDG-2024,Boroczky-Lutwak-Yang-Zhang-2013}. Here $\lambda_K$ denotes the  principal curvature radii of $\partial K$, $E_k(\lambda_K)=\frac{\sigma_k(\lambda_K)}{\binom{n}{k}}$ (see \eqref{equ:sigma_k}) is the normalized $k$-th elementary symmetric polynomial of $\lambda_K$, and $dx:=d\mu_{\mathbb{S}^n}(x)$ is the standard induced volume element on $\mathbb{S}^n$.

The behavior below $p=1$ changes sharply with the parameter. For $0<p<1$, as shown by counterexamples in \cite{Boroczky-Lutwak-Yang-Zhang-2012}, the $L^p$-Brunn--Minkowski inequality and the $L^p$-Minkowski inequality do not hold in general. However, it was conjectured and later proved partially in \cite{Boroczky-Lutwak-Yang-Zhang-2012,Kolesnikov-Milman-MAMS-2022,Chen-Huang-Li-Liu-2020,Putterman-2021} that it holds under the additional assumption that $K$ and $L$ are origin-symmetric.
The case $p=0$ in \eqref{equ:BM-neq-p=1} is interpreted in the limiting sense as $p\to 0$ and is known as the log-Brunn--Minkowski inequality
\begin{align}\label{equ:BM-neq-p=0}
	V_n((1-t)\cdot K+_0t\cdot L)\geq V_n(K)^{1-t}V_n(L)^{t}.
\end{align} This inequality was proved for origin-symmetric convex bodies in dimension two in \cite{Boroczky-Lutwak-Yang-Zhang-2012}. In higher dimensions and in the non-origin-symmetric setting, partial results under additional assumptions were obtained in \cite{Xi-Leng-2016,Boroczky-Kalantzopoulos-2022}.


For $p<0$, the situation is more difficult still, even within the origin-symmetric class. The counterexamples in \cite{Boroczky-Saroglou-Cvpde-2024,Kolesnikov-Milman-MAMS-2022,Milman-JDG-2024} showed that the $L^p$-Brunn--Minkowski and $L^p$-Minkowski inequalities fail for some origin-symmetric bodies. Moreover, \cite{Kolesnikov-Milman-MAMS-2022,Boroczky-Saroglou-Cvpde-2024} established nonuniqueness for the even $L^p$-Minkowski problem for certain $p<0$. Here the $L^p$-Minkowski problem asks for a convex body whose $L^p$-surface area measure is a prescribed measure, while ``even'' means that 
the solutions are restricted to origin-symmetric convex bodies.
More precisely, \cite{Boroczky-Saroglou-Cvpde-2024} proved that, when $-(n+1)<p<0$, the set of nonunique solutions to the $L^p$-Minkowski problem is dense in the class of convex bodies containing the origin. Consequently, for $p<0$, the $L^p$-Brunn--Minkowski inequality fails for generic origin-symmetric convex bodies.





The $L^p$-Brunn--Minkowski inequality has inspired many generalizations, often referred to as Brunn--Minkowski-type inequalities, in which volume is replaced by other geometric functionals. For instance, partial results have been obtained for the dual quermassintegral by \cite{Sadovsky-ZhangGaoYong-2025,Xi-Zhang-2022}, and for the intrinsic volume $V_k$ (defined in \eqref{equ:intrinsic volume}) by \cite{Bianchini-Colesanti-Pagnini-Roncoroni-Math-ann-2023, Patsalos-Saroglou-arxiv}. On the other hand, 
\cite{Patsalos-Saroglou-arxiv}*{Theorem 1.1}
provided counterexamples showing that, for $p\in[0,1)$ and every $k\in\{1,\dots,n-1\}$, the corresponding $L^p$-Brunn--Minkowski inequality for $V_k$ may fail even for origin-symmetric convex bodies.

For origin-symmetric convex bodies $K,L$ and $p>-(n+1)$, the equivalence among the $L^p$-Brunn--Minkowski inequality, uniqueness for the $L^p$-Minkowski problem, and the $L^p$-Minkowski inequality was proved in \cite{Boroczky-Lutwak-Yang-Zhang-2012,Milman-JEMS-2025}. Thus
in light of these counterexamples, establishing the $L^p$-Brunn--Minkowski inequality is highly challenging and likely impossible in the following cases: for $0\leq p<1$ without origin-symmetry, for $p < 0$ even under origin-symmetry, and in the setting of intrinsic volumes. In this paper, somewhat surprisingly, we provide affirmative answers in all of the above cases and further establish strengthened versions of the corresponding inequalities, namely, their stability.

As noted above, even within the class of origin-symmetric convex bodies, these inequalities still fail in the vast majority of cases. The natural way to recover positive results is therefore to impose a local geometric hypothesis. Thus,
inspired by \cite{Boroczky-Chen-Liu-Saroglou-arxiv,Boroczky-Saroglou-Cvpde-2024,Patsalos-Saroglou-arxiv,Bianchini-Colesanti-Pagnini-Roncoroni-Math-ann-2023}, we restrict attention to convex bodies that are $C^2$-close to the unit ball. More precisely, let $\mathcal{N}_{\varepsilon}$ be the collection of all $C^2$ convex bodies $K$ satisfying
\begin{align}\label{equ:def-Ne}
	\| u_K - 1 \|_{C^2(\mathbb{S}^n)} < \varepsilon, \quad \text{for some}\ \varepsilon>0.
\end{align}
Here we use the notation
$$\| f \|_{C^m}:=\| f \|_{C^m(\mathbb{S}^n)}=\max_{0\leq k\leq m}\sup_{x\in\mathbb{S}^n}\left|\nabla^k f(x)\right|$$
for $m\in\mathbb{Z}^+$,
where $\nabla$ denotes the Levi-Civita connection of the induced metric on $\mathbb{S}^n$.

Under this local condition, Patsalos and Saroglou \cite{Patsalos-Saroglou-arxiv} proved an $L^p$-Brunn--Minkowski inequality for intrinsic volumes of origin-symmetric convex bodies. Their theorem strengthens the one in \cite{Bianchini-Colesanti-Pagnini-Roncoroni-Math-ann-2023} not only by replacing the condition that  $L$ is a ball with the weaker requirement that $L\in \mathcal{N}_{\varepsilon}$ is an origin-symmetric convex body, but also by allowing $p$ to take negative values within a certain interval.
\begin{theorem}[Patsalos--Saroglou, \cite{Patsalos-Saroglou-arxiv}*{Theorem 1.2}]
	\label{thm:Arxiv-origin-symmeric-BM}
	Fix $n\geq 2$, $k\in\{1,\cdots,n\}$, and $p\in (-k+1-\frac{2k}{n},1)$. There  exists $\varepsilon>0$ that depends only on $n,p$, such that for all $t\in(0,1)$ and for every \textbf{origin-symmetric} convex bodies $K,L\in\mathcal{N}_{\varepsilon}$, it holds
	\begin{align}
		\label{equ:BMneq-V_k}
		V_k\left((1-t) \cdot K+_p t \cdot L\right) \geq\left((1-t) V_k(K)^{\frac{p}{k+1}}+t V_k(L)^{\frac{p}{k+1}}\right)^{\frac{k+1}{p}}.
	\end{align}
	Equality holds if and only if $K$ and $L$ are dilates.
\end{theorem}

Let us introduce the intrinsic volume $V_k(\cdot)$ appearing in Theorem  \ref{thm:Arxiv-origin-symmeric-BM}.
For convex bodies $K, L \subset \mathbb{R}^{n+1}$, the mixed volume (cf. \cite{Xia2017, book-convex-body}) is defined  by $$V_{k}(K,L)=V(\underbrace{K, \cdots, K,}_{(n-k+1) \text { copies }} \underbrace{L, \cdots, L}_{k \text { copies }}),$$ with  $k=0,\cdots,n+1$. In particular, $V_{n+1}(K,L)=(n+1)|L|$ and $V_0(K,L)=(n+1)|K|$. The \textbf{$k$-th intrinsic volume} is
\begin{align}\label{equ:intrinsic volume}
	V_k(K)=V_{n-k}(K,\mathbb{B}^{n+1}),
\end{align}
where $-1\leq k\leq n$. Specially, when $k=n$, $V_n(K)=(n+1)|K|$.
If $K$ is of class $C^2$, then $V_k(K)$ can be expressed in the following form:
\begin{align}
	V_k(K)=\int_{\mathbb{S}^n}u_KE_k(\lambda_K)~dx=\int_{\mathbb{S}^n}E_{k+1}(\lambda_K)~dx.
\end{align}

In this paper, we pursue two directions of extension. First, we consider the non-origin-symmetric setting. The origin-symmetry assumption is essential in \cite{Patsalos-Saroglou-arxiv}, as the authors need to determine the eigenvalues of an operator on the class of origin-symmetric convex bodies. Consequently, the approach adopted in this paper is fundamentally different from that in \cite{Patsalos-Saroglou-arxiv}. Second, we establish strengthened versions of the inequalities, namely, their stability. Thus, rather than merely proving an $L^p$-Brunn--Minkowski inequality for intrinsic volumes, we measure how close $K$ and $L$ must be to homothetic bodies when equality is nearly attained. To quantify this closeness,
we introduce the following two distance measures to quantify the deviation between $K$ and $L$:
(i) the classical
$L^2$ distance (cf. \cite{book-convex-body}*{eq. (7.120)})
\begin{align}\label{equ:L2-distance}
	\delta_2(K,L)=\left\Vert u_K-u_L\right\Vert_{L^2}=\left(\int_{\mathbb{S}^n}|u_K-u_L|^2~dx\right)^{\frac{1}{2}},
\end{align}
and notably,
(ii) the more natural \textbf{relative asymmetry index} used by Figalli, Maggi, and Pratelli \cite{Figalli-Maggi-Pratelli-2010}*{eq.  (1.29)}, which is defined by
\begin{align}
	\label{equ:def-asymmetry}
	\alpha(K,L)=\min_{y \in \mathbb{R}^{n+1}}\left\{\frac{|K \Delta(r L+y)|}{|K|}:|r L|=|K|\right\},
\end{align}
where  $K \Delta L=(K \backslash L) \cup(L \backslash K)$ denotes the symmetric difference between two sets $K,L\subset\mathbb{R}^{n+1}$.
Distance $\alpha(K,L)$ satisfies symmetry and the triangle inequality \cite{Figalli-Maggi-Pratelli-2010}*{eq. (4.2)}.
For the case $u_L\equiv1$, i.e., $L=\mathbb{B}^{n+1}$, the relative asymmetry index reduces to the Fraenkel asymmetry (ref. \cite{Fusco-Maggi-2008}).

If $k=n$, the stability of $L^p$-Brunn--Minkowski inequalities is a well-known problem that has attracted considerable attention in recent years. When $p=1$, it reduces to the stability of the classical Brunn--Minkowski inequality.
The first stability forms of the Brunn--Minkowski inequality for convex
bodies were due to Minkowski himself  and more details can be seen in \cite{Groemer-1993}. Subsequently, nearly optimal stability results were obtained by Diskant \cite{Diskant-1973} and Groemer \cite{Groemer-1988}. Later developments introduced natural distance notions and established the optimal exponent, with further extensions to general bounded measurable sets in \cite{Figalli-Maggi-Pratelli-2010,Figalli-Maggi-Pratelli-2009,Figalli-Hintum-Tiba-arxiv-2023,Figalli-Hintum-Tiba-arxiv-2024}.
For $p>1$, B\"{o}r\"{o}czky, Machado, and Ramos \cite[Lemma 4.1]{Boroczky-Machado-Ramos-2026-arxiv} established stability of the $L^p$-Brunn--Minkowski inequality by invoking the quantitative $p=1$ result  from \cite{Boroczky-Figalli-Ramos-2025-book}. 
At $p=0$, B\"{o}r\"{o}czky and De \cite{Boroczky-De-2024} proved stability under symmetry with respect to $n+1$ independent linear hyperplanes. In \cite{Colesanti-Livshyts-Livshyts-2017}, Colesanti, Livshyts, and Marsiglietti established the stability results near a Euclidean ball for the dimensional Brunn--Minkowski inequality for radially symmetric log-concave measures and for the log-Brunn--Minkowski inequality.
By contrast, to the best of our knowledge, no stability result was previously available for $0\leq p<1$ without origin symmetry, for $p<0$ even under origin symmetry, or for intrinsic volumes. Somewhat surprisingly, this paper provides affirmative results in all three difficult regimes identified above and yields their stronger stability forms.

We consider a more general setting where the convex bodies
$K,L$ are not required to be origin-symmetric. An additional difficulty in the non-origin-symmetric setting is the choice of a centering condition compatible with $L^p$ interpolation. For this purpose, we introduce the \textbf{$p$-weighted Steiner point} of a convex body $K$:
\begin{equation}
	\label{equ:def-p-Steiner}
	\mathcal{S}_p(K)=
	\left\lbrace
	\begin{array}{ll}
		\int_{\mathbb{S}^n}u_K^p(x)\cdot x~dx, &\text{\ for\ } p\neq 0,
		\\[3pt]
		
		\int_{\mathbb{S}^n}\log u_K(x)\cdot x~dx, &\text{\ for\ } p= 0.
	\end{array}
	\right.
\end{equation}

Throughout the paper, we impose the normalization $\mathcal{S}_p(K)=\mathcal{S}_p(L)$, and  assume $\mathcal{S}_p(K)=\mathcal{S}_p(L)=o$ in Section \ref{sec 5}. The following observations explain both its scope and its compatibility with $L^p$-Minkowski interpolation.

 The condition $\mathcal{S}_p(K)=\mathcal{S}_p(L)$ is necessary. Indeed, if we take $K=\mathbb{B}^{n+1}+\frac{\varepsilon}{2} e_1$ and $L=\mathbb{B}^{n+1}-\frac{\varepsilon}{2} e_1$, where $e_1=(0,\cdots,0,1)\in\mathbb{R}^{n+1}$, then the $L^p$-Minkowski and $L^p$-Brunn--Minkowski inequalities for intrinsic volumes and volume fail for such $K,L\in \mathcal{N}_{\varepsilon}$ and $0<p<1$.
\begin{remark}
	(1) When $p=1$, \eqref{equ:def-p-Steiner} reduces to the classical Steiner point \cite{book-convex-body}, denoted by $\mathcal{S}(K):=\mathcal{S}_1(K)$. When $p=-n-2$, it reduces to the definition of the Santal\'o point.
	
	(2) The condition $\mathcal{S}_p(K)=o$ contains the origin-symmetric case as a special case.
	
	(3) The definition of $\mathcal{S}_p(K)$ follows the weighted centering conditions for the entropy point in \cite{Andrews-Guan-Ni-2016}*{Lemma 2.5, eq. (2.7)} and \cite{Guan-Ni-2017}*{Lemma 2.2, eq. (2.2)}. It is this normalization that allows us to treat the general case without assuming origin symmetry.

	
	(4) Even if the Steiner points of $K,L\in\mathcal{N}_{\varepsilon}$ lie at the origin, this does not in general imply that the Steiner point of their $L^p$-Minkowski convex combination
	$Q_t:=(1-t)\cdot K+_pt\cdot L\in\mathcal{N}_{O(\varepsilon)}$ with $t\in[0,1]$ is also at the origin, except in the special case $p=1$.
	In contrast, a key feature of $\mathcal{S}_p(K)=\mathcal{S}_p(L)=o$ is that it is preserved under $L^p$-Minkowski convex combinations. That is, if $K,L\in\mathcal{N}_{\varepsilon}$  satisfy $\mathcal{S}_p(K)=\mathcal{S}_p(L)=o$, then $\mathcal{S}_p(Q_t)=o$. This follows directly from Lemma \ref{lemma:C2-1-p-sum-C2-1}.
	

	(5) The condition $\mathcal{S}_p(K)=\mathcal{S}_p(L)$ is not preserved under translations of $K$ or $L$. This is consistent with the fact that translations do not appear in the equality case of the $L^p$-Minkowski inequality stated in Corollary \ref{thm:Christoffel} (a). 
\end{remark}



Under this centering condition, we establish stability of the $L^p$-Brunn--Minkowski inequality for non-origin-symmetric convex bodies over a wider range of $p$, including negative values, and with respect to several notions of distance. These results extend \cite{Bianchini-Colesanti-Pagnini-Roncoroni-Math-ann-2023} and Theorem \ref{thm:Arxiv-origin-symmeric-BM}, while relying on a method different from those used in the references discussed above. Our first main theorem, derived from Corollary  \ref{thm:p-M-L2}, gives $L^2$-stability simultaneously for volume and intrinsic volumes.

\begin{theorem}
	\label{thm:p-BM-L2-Vk=Vk}
	Let $K, L \in \mathcal{K}^{n+1}_o$ be $C^3$  convex bodies with $V_k(K)=V_k(L)$, $1\leq k\leq n$. For every $\eta>0$, there exists $\varepsilon>0$  such that if $K,L\in\mathcal{N}_{\varepsilon}$ and $\mathcal{S}_p(K)=\mathcal{S}_p(L)=o$, then for $\forall t\in[0,1]$ we have for $p=0$,
	\begin{align}\label{equ:p-BM-Vk=p=0}
		&\log	\left(V_k\((1-t)\cdot K+_0t\cdot L\)\right)-\log\left(V_k(K)\right)\nonumber
		\\
		\geq&
		(k+1)t(1-t)V_k(K)^{-1}\left(\frac{1}{2}(k-1+\frac{2k}{n})-\eta\right)\delta_2^2(K,L);
	\end{align}
	and for $p\in (-k+1-\frac{2k}{n},0)\cup(0,+\infty)$,
	\begin{align}\label{equ:p-BM-Vk=p!=0}
		&\frac{1}{p}\left(V_k\((1-t)\cdot K+_pt\cdot L\)^{\frac{p}{k+1}}-V_k(K)^{\frac{p}{k+1}}\right)\nonumber
		\\
		\geq&t(1-t)V_k(K)^{-\frac{k+1-p}{k+1}}\left(\frac{1}{2}(p+k-1+\frac{2k}{n})-\eta\right)\delta_2^2(K,L).
	\end{align}
\end{theorem}
\begin{remark}\label{rk:L2-replaced by the anisotropic weighted-BM}
	Theorem \ref{thm:p-BM-L2-Vk!=Vk} gives the corresponding $L^2$-stability estimate for the unnormalized $L^p$-Brunn--Minkowski inequality. Here ``unnormalized'' means that $V_k(K)=V_k(L)$ is not assumed.
\end{remark}
In terms of the relative asymmetry index, Theorem \ref{thm:p-BM-L2-Vk=Vk-index} provides stability results for intrinsic volumes and  volume, which are derived  from Corollary \ref{thm:p-M-index-k!=n=!Vk}.
\begin{theorem}
	\label{thm:p-BM-L2-Vk=Vk-index}
	Let $K, L \in \mathcal{K}^{n+1}_o$ be  $C^3$  convex bodies with $V_k(K)=V_k(L)$, $1\leq k\le n$. For all $\eta>0$, there exists $\varepsilon>0$  such that if $K,L\in\mathcal{N}_{\varepsilon}$ and $\mathcal{S}_p(K)=\mathcal{S}_p(L)=o$, then for $\forall t\in[0,1]$ we have for $p=0$,
	\begin{align}\label{equ:p-BM-Vk=p=0-index}
		&\log	\left(V_k\((1-t)\cdot K+_0t\cdot L\)\right)-\log\left(V_k(K)\right)\nonumber
		\\
		\geq&
		(k+1)t(1-t)\left(\frac{(k-1)n+2k}{2n(n+1)^2}-\eta\right)\alpha(K,L)^2;
	\end{align}
	and for $p\in (-k+1-\frac{2k}{n},0)\cup(0,+\infty)$,
	\begin{align}\label{equ:p-BM-Vk=p!=0-index}
		&\frac{1}{p}\left(V_k\((1-t)\cdot K+_pt\cdot L\)^{\frac{p}{k+1}}-V_k(K)^{\frac{p}{k+1}}\right)\nonumber
		\\
		\geq&t(1-t)V_k(K)^{\frac{p}{k+1}}\left(\frac{n(k-1+p)+2k}{2n(n+1)^2}-{\eta}\right)\alpha(K,L)^2.
	\end{align}
\end{theorem}
\begin{remark}\label{rk:L2-replaced by the anisotropic weighted-BM-index}
	(1) Theorem \ref{thm:p-BM-L2-Vk!=Vk-index} gives stability in terms of the relative asymmetry index for the unnormalized $L^p$-Brunn--Minkowski inequality, that is, without requiring $V_k(K)=V_k(L)$.
	
	(2) For $p>1$ and $k=n$, B\"{o}r\"{o}czky, Machado, and Ramos \cite[Lemma 4.1]{Boroczky-Machado-Ramos-2026-arxiv} obtained stability for arbitrary convex bodies containing the origin by using the $p=1$ stability result. Their coefficient is $\theta_n=cn^{-4}(\log n)^{-1}$, with $c\in(0,1)$. Theorem \ref{thm:p-BM-L2-Vk=Vk-index}  improves this coefficient under the additional assumption that $K$ and $L$ are sufficiently close to a ball.
	
	(3) The stability results from Theorem \ref{thm:p-BM-L2-Vk=Vk} (or Theorem \ref{thm:p-BM-L2-Vk=Vk-index}) implies the $L^p$-Brunn--Minkowski inequality for $p > -k+1-\frac{2k}{n}$, since $\delta_2^2(K,L)$ and $\alpha(K,L)^2$ are non-negative. In other words, stability contains the original inequality as a special case. 
\end{remark}





The stability of  $L^p$-Minkowski inequality for volume and intrinsic volumes (see Corollaries \ref{thm:p-M-L2} and \ref{thm:p-M-index-k!=n=!Vk} for details)
can imply the corresponding $L^p$-Brunn--Minkowski inequality. 
Accordingly, to obtain stability for the $L^p$-Brunn--Minkowski inequality, we first establish the  corresponding stability results for the $L^p$-Minkowski inequality in Section \ref{sec 4}.


Since $\delta_2(N_1K, N_2L) \ge 0$, with equality if and only if $K = cL$ for some constant $c>0$, the stability results in Corollary \ref{thm:p-M-L2} and Theorem \ref{thm:p-BM-L2-Vk!=Vk} yield the following $L^p$-Minkowski inequality (Corollary \ref{thm:Christoffel} (a)) and $L^p$-Brunn--Minkowski inequality (Corollary \ref{thm:Christoffel} (b)) for intrinsic volumes and volume, valid for $p > -k+1-\frac{2k}{n}$. Moreover, the $L^p$-Minkowski inequality for intrinsic volumes can be used to derive the uniqueness for the $L^p$-Christoffel-Minkowski problem as follows (Corollary \ref{thm:Christoffel} (c)) .
\begin{corollary}\label{thm:Christoffel}
	Let $1\leq k\leq  n$ and $p>-k+1-\frac{2k}{n}$. There exists $\varepsilon>0$, that depends only on $k,\, n$ and $p$, such that
	if $K,L\in\mathcal{N}_{\varepsilon}$ are $C^3$ convex bodies satisfying $\mathcal{S}_p(K)=\mathcal{S}_p(L)=o$,
	then
	
	(a) we have the $L^p$-Minkowski inequality,
	\begin{align}
		\label{equ:p-M-inequ-p=0}
		\int_{\mathbb{S}^n}\left(\log(\frac{u_L}{u_K})\right)u_KE_k(\lambda_K) ~dx
		\geq
		\frac{1}{k+1}V_k(K)\log\frac{V_k(L)}{V_k(K)}, \quad \mathrm{ if }\ p=0,
	\end{align}
	\begin{align}
		\label{equ:p-M-inequ-p!=0}
		\frac{1}{p}\int_{\mathbb{S}^n}u_L^pu_K^{1-p}E_k(\lambda_K)~dx\geq \frac{1}{p}V_k(K)^{\frac{k+1-p}{k+1}}V_k(L)^{\frac{p}{k+1}}, \quad \mathrm{ if }\ p\ne0;
	\end{align}
	
	(b) we have the $L^p$-Brunn--Minkowski inequality,
	\begin{align}
		\label{equ:p-BM-inequ-p=0}
		V_k\((1-t)\cdot K+_0t\cdot L\)
		\geq  V_k(K)^{1-t} V
		_k(L)^{t}, \quad \mathrm{ if }\ p=0,
	\end{align}
	\begin{align}
		\label{equ:p-BM-inequ-p!=0}
		V_k\((1-t)\cdot K+_pt\cdot L\)\geq \left((1-t)V_k(K)^{\frac{p}{k+1}}+tV_k(L)^{\frac{p}{k+1}}\right)^{\frac{k+1}{p}}, \quad \mathrm{ if }\ p\ne0;
	\end{align}
	
	(c) and the $L^p$-Christoffel-Minkowski equation	
	\begin{align}\label{equ:Christolell-equ}
		u_{K}^{1-p} E_k\left(\lambda_K\right)=u_{L}^{1-p} E_k\left(\lambda_L\right)
	\end{align}
	implies $K=L$ when $p\neq 1+k$ and   $K=NL$ when $p= 1+k$.
	
	Moreover,  equality holds in \eqref{equ:p-M-inequ-p=0}
	$\thicksim$\eqref{equ:p-BM-inequ-p!=0}, if and only if $K=cL$ for some constant $c>0$. 
	
\end{corollary}


\begin{remark}
	 The case $k=n$ in Corollary \ref{thm:Christoffel} (c) corresponds to the uniqueness of $L^p$-Minkowski problem, while the case $1\leq k<n$ corresponds to the uniqueness of  $L^p$-Christoffel-Minkowski problem.
\end{remark}

Li and Wan \cite{Li-Wan-2024} and  Guan and Xia \cite{Xia-Guan-2018} established uniqueness results for the $L^p$-Christoffel--Minkowski problem when $p>1$. The uniqueness problem for the $L^p$-Christoffel--Minkowski problem when $p<1$ is considerably more difficult. Under the assumptions that $p\in (-k+1-\frac{2k}{n},1)$, the prescribed function $f$ is even, and
$
\|f-1\|_{C^\alpha}<\delta,
$
uniqueness among origin-symmetric convex bodies was proved in \cite{Patsalos-Saroglou-arxiv}*{Theorem 1.4}.
In the more special case $f\equiv 1$, uniqueness without assuming origin symmetry was established in \cite{Chen-JFA-2020}*{Theorem 1.1} for $1-k<p<1$ and in \cite{McCoy-2011} for $p=1-k$.

The rest of this paper is organized as follows. In Section \ref{sec 2}, we briefly introduce some preliminaries of mixed discriminants.
In Section \ref{sec 3}, we study  some characteristics and formulations for  convex bodies that are $C^2$-close to a ball.
In Section \ref{sec 4}, we   establish the stability of $L^p$-Minkowski inequality, to prove Theorems \ref{thm:p-M-L2-Vk=Vk} and  \ref{thm:p-M-index-kneqn}, as well as Corollaries \ref{thm:p-M-L2} and \ref{thm:p-M-index-k!=n=!Vk}.
Section \ref{sec 5} concerns the stability of the $L^p$-Brunn--Minkowski inequality (see Theorems \ref{thm:p-BM-L2-Vk=Vk},  \ref{thm:p-BM-L2-Vk=Vk-index},  \ref{thm:p-BM-L2-Vk!=Vk}, and \ref{thm:p-BM-L2-Vk!=Vk-index}), which follows from Corollaries \ref{thm:p-M-L2} and \ref{thm:p-M-index-k!=n=!Vk}. Finally, in Section \ref{sec 6}, we obtain a uniqueness result for the $L^p$-Christoffel-Minkowski problem (see Corollary \ref{thm:Christoffel} (c)).

\section{Preliminary}\label{sec 2}
In this section, we introduce the mixed discriminants.
\begin{definition}[\cite{book-convex-body}*{Section 5.5}, \cite{Xia-arxiv}*{Page 10}, and \cite{Ivaki-Milman-2023}*{Page 4}]
	The mixed discriminant $\mathcal{Q}:\left(\mathbb{R}^{n \times n}\right)^n \rightarrow \mathbb{R}$ is defined by	
	\begin{align}
		\operatorname{det}\left(\lambda_1 A_1+\cdots+\lambda_m A_m\right)=\sum_{i_1, \cdots, i_n=1}^m \lambda_{i_1} \cdots \lambda_{i_n} \mathcal{Q}\left(A_{i_1}, \cdots, A_{i_n}\right),\label{equ:QA}
	\end{align}
	for $m \in \mathbb{N}, \lambda_1, \cdots, \lambda_m \geq 0$ and the real symmetric matrices $A_1, \cdots, A_m \in \mathbb{R}^{n \times n}$. If $\left(A_k\right)_{i j}$ denotes the $(i, j)$-element of the matrix $A_k$, then	
	\begin{align}\label{equ:Q1...n}
		\mathcal{Q}\left(A_1, \cdots, A_n\right)=\frac{1}{n!} \sum_{i_1, \cdots, i_n; j_1, \cdots, j_n} \delta_{j_1 \cdots j_n}^{i_1 \cdots i_n}\left(A_1\right)_{i_1 j_1} \cdots\left(A_n\right)_{i_n j_n}.
	\end{align}	
\end{definition}
From \eqref{equ:QA}, it is easy to check that for any invertible matrix $B\in \mathbb{R}^{n \times n}$, there holds
\begin{align}\label{equ:Q(AB)=Q(A)det(B)}
	\mathcal{Q}\left(A_1B, \cdots, A_nB\right)=\mathcal{Q}\left(A_1, \cdots, A_n\right)\det(B).
\end{align}

From \eqref{equ:Q1...n}, we check that $\mathcal{Q}(A_1,\cdots,A_n)$ is completely symmetric in the variables $A_1,\cdots,A_n$, and  for any $a,b\in\mathbb{R}$, we have
\begin{align}\label{equ:Q(A+B)}
	\mathcal{Q}\left(a\cdot A+b\cdot B, A_2,\cdots, A_n\right)=a\cdot\mathcal{Q}\left(A, A_2,\cdots, A_n\right)+b\cdot\mathcal{Q}\left(B, A_2,\cdots, A_n\right),
\end{align}
and for $A_1\neq A_2,\cdots,A_n$, we have
\begin{align}
	\label{equ:Q=Qij*Aij}
	\mathcal{Q}(A_1,\cdots,A_n)=\sum_{i,j=1}^{n}(A_1)_{ij}\cdot \frac{\partial}{\partial (A_1)_{ij}} \mathcal{Q}\left(A_1, \cdots, A_n\right).
\end{align}

For convenience, we denote
\begin{align}
	&\mathcal{Q}\left(A_1[k_1],\cdots, A_{r-1}[k_{r-1}],A_r\right)
	\nonumber
	\\
	:=&\mathcal{Q}\left(A_1[k_1],\cdots, A_{r-1}[k_{r-1}],A_r[n-k_1-\cdots-k_{r-1}]\right)
	\nonumber
	\\
	:=&\mathcal{Q}\(\underbrace{A_1\cdots,A_1}_{k_1\text{\ copies}},\cdots, \underbrace{A_{r-1},\cdots,A_{r-1}}_{k_{r-1}\text{\ copies}}, \underbrace{A_{r},\cdots,A_{r}}_{(n-k_1-\cdots-k_{r-1})\text{\ copies}}\).
	\label{equ:Q(A[])}
\end{align}

A useful identity of elementary symmetric function $\sigma_k$ is
\begin{align}
\sigma_k(\lambda(A))=\sigma_k(A):=&\frac{1}{k!} \sum_{i_1, \cdots, i_k; j_1, \cdots, j_k} \delta_{j_1 \cdots j_k}^{i_1 \cdots i_k}A_{i_1 j_1} \cdots A_{i_k j_k}\label{equ:sigma_k}
	\\
	=&\binom{n}{k}\cdot \mathcal{Q}(A[k],I[n-k]),\label{equ:sigma_k=cnkQ}
\end{align}
where  $\binom{n}{k}=\frac{n!}{k!(n-k)!}$, and $\lambda(A)$ denotes the eigenvalue vector of $A$.
Similarly, it is easy to check
\begin{align*}
	\frac{\partial}{\partial (A_1)_{ij}}\mathcal{Q}\left(A_1, \cdots, A_n\right)=\frac{1}{n!} \sum_{i_2, \cdots, i_n;  j_2, \cdots, j_n} \delta_{jj_2 \cdots j_n}^{ii_2 \cdots i_n}\left(A_2\right)_{i_2 j_2} \cdots\left(A_n\right)_{i_n j_n},
\end{align*}
where $A_1\neq A_2,\cdots,A_n$.
Then for $A\neq A_{k+1},\cdots,A_n$, we have 
\begin{align}
	&\frac{\partial}{\partial (A)_{ij}}\mathcal{Q}\left(A[k], A_{k+1}\cdots, A_n\right)
	\nonumber
	\\
	=&\frac{k}{n!} \sum_{ i_2, \cdots, i_n;  j_2, \cdots, j_n} \delta_{jj_2 \cdots j_n}^{ii_2 \cdots i_n}A_{i_2j_2}\cdots A_{i_kj_k}\left(A_{k+1}\right)_{i_{k+1} j_{k+1}} \cdots\left(A_n\right)_{i_n j_n}\label{equ:pQ1...n}
\end{align}
We denote by $\tau[f]=\{\tau_{ij}[f]\}_{n\times n}$ the matrix with entries $$\tau_{ij}[f]:=\nabla_{e_i}\nabla_{e_j}f+f\delta_{ij}$$ for any $f\in C^3(\mathbb{S}^n)$, where  $\{e_i\}_{i=1}^n$ is a local orthonormal frame on $\mathbb{S}^n$ and $\nabla$ denotes the Levi-Civita connection of the induced metric on $\mathbb{S}^n$.  By \cite{Andrews-97}*{Proposition 3},  for any $f_i\in C^3(\mathbb{S}^n)$ with $i=1,\cdots, n+1$, we have
\begin{align}
	\label{equ:Q_ij,j=0}
	\sum_{j=1}^{n}\nabla_j\frac{\partial\mathcal{Q}\left(\tau[f_1],\cdots,\tau[f_n]\right)}{\partial \tau_{ij}[f_1]}=0,
\end{align}
and if $\tau[f_i]>0$ with $i=1, \cdots, n+1$, then
\begin{align}
	\label{equ:int-uQ-sysmetry}
	\int_{\mathbb{S}}f_1\mathcal{Q}\left(\tau[f_2],\tau[f_3],\cdots, \tau[f_{n+1}]\right)~dx=\int_{\mathbb{S}}f_2\mathcal{Q}\left(\tau[f_1],\tau[f_3],\cdots, \tau[f_{n+1}]\right)~dx.
\end{align}
 By the linearity of $\mathcal{Q}$ with respect to each $f_i$ (for all $i\in\{1,\cdots,n+1\}$), it is easy to see that, in the above formula \eqref{equ:int-uQ-sysmetry}, $f_i$ with positive definite $\tau[f_i]$ can also be replaced by $f-f'$, where both $\tau[f]$ and $\tau[f']$ are positive definite.
There is also a more general version of the definitions and properties discussed above. The capillary boundary case is covered in \cite{Xia-arxiv}, and the anisotropic capillary boundary case is covered in \cite{Gao-Li-cvpde}.


In this paper, we consider convex bodies $K, L \in \mathcal{N}_{\varepsilon}$ (i.e., satisfying \eqref{equ:def-Ne}). Their support functions are denoted by $u_K = 1 + w$ and $u_L = 1 + \bar{w}$, respectively. For sufficiently small $\varepsilon$, we have $\tau[u_K] > 0$ and $\tau[u_L] > 0$, which ensures that $K$ and $L$ are indeed convex bodies.


\section{Some propositions for convex bodies that are $C^2$-close to a ball}\label{sec 3}

In this section, we consider  convex bodies that are $C^2$-close to a ball, and compute the corresponding geometric quantities.


We reformulate the  curvature functions $E_k(\lambda_K)$ and intrinsic volumes $V_k(K)$ of $K$ as follows:

\begin{lemma}
	\label{lemma:I_k}
	For convex body $ K$ with support function $u_K(x)=1+w(x)\in C^2(\mathbb{S}^n)$, we have
	\begin{align}
		\label{equ:E_k(lambda)}         E_k(\lambda_K)=\sum_{s=
			0}^{k}\binom{k}{s}\mathcal{Q}\(\nabla^2 w[s],I\)(1+w)^{k-s};
	\end{align}
	and
	\begin{align}
		\label{equ:Vk}
		V_k(K)=\int_{\mathbb{S}^n} \sum_{s=
			0}^{k}\binom{k}{s}\mathcal{Q}\(\nabla^2 w[s],I\)(1+w)^{k-s+1}~dx.
	\end{align}
	Here $k=0,\cdots,n$ and the mixed discriminant  $\mathcal{Q}$ is defined by \eqref{equ:Q1...n} and \eqref{equ:Q(A[])}.
\end{lemma}
\begin{proof}
	Let $\lambda_K$ denote the  principal curvature radii  for hypersurface
	$\partial K$. It is well known (see, e.g., \cite[eq. (2.1)]{Ding-Li-2025-JDE}) that $\lambda_K$ is given by the eigenvalues of the matrix
	\begin{align}\nonumber
		\tau[u_K]=\{\nabla^2u_K+u_KI\}=\{\nabla^2w+(1+w)I\}.
	\end{align}
	Combining  above with  
	\eqref{equ:sigma_k=cnkQ} and \eqref{equ:Q(A+B)}, we obtain
	\begin{align*}
		E_k(\lambda_K)=&\frac{\sigma_k(\lambda_K)}{\binom{n}{k}}\overset{\eqref{equ:sigma_k=cnkQ}}{=}\mathcal{Q}\(\left(\nabla^2w+(1+w)I\right)[k],I\)
		\\
		\overset{\eqref{equ:Q(A+B)}}{=}&\sum_{s=
			0}^{k}\binom{k}{s}\mathcal{Q}\(\nabla^2 w[s],(1+w)I[k-s],I\)
		\\
		=&\sum_{s=
			0}^{k}\binom{k}{s}\mathcal{Q}\(\nabla^2 w[s],I\)(1+w)^{k-s},
	\end{align*}
which is \eqref{equ:E_k(lambda)}. By definition,
	\begin{align*}
		V_k(K)=&\int_{\mathbb{S}^n}u_KE_k(\lambda_K)~dx
		\\
		=&\int_{\mathbb{S}^n} \sum_{s=
			0}^{k}\binom{k}{s}\mathcal{Q}\(\nabla^2 w[s],I\)(1+w)^{k-s+1}~dx,
	\end{align*}
and then	\eqref{equ:Vk} holds.
\end{proof}

 For convenience, we write $$v:=u_L-u_K,$$ and denote $$\mathscr{E}:=O(\varepsilon)\cdot \(\Vert v\Vert_{L^2}^2+\Vert \nabla v\Vert_{L^2}^2\).$$
	\begin{lemma}
		\label{lemma:int vEk<OE}
		Let $1\leq k\leq n$, and $K,L\in\mathcal{N}_{\varepsilon}$ be  $C^3$  convex bodies with $V_k(K)=V_k(L)$. Then we have 
		\begin{align}
			\label{equ:int vEk<OE}
			\int_{\mathbb{S}^n}vE_k(\lambda_K)~dx=-\frac{k}{2}\Vert v\Vert_{L^2}^2+\frac{k}{2n}\Vert \nabla v\Vert_{L^2}^2+\mathscr{E}.
		\end{align}
	\end{lemma}
	\begin{proof}
		For any $t\in[0,1]$, we take $u_t=u_K+tv=tu_L+(1-t)u_K$, then $u_0=u_K$, $u_1=u_L$,  $\frac{d}{dt}u_t=v$ and $\frac{d}{dt}\tau[u_t]=\tau[v]$. Since $K,L\in\mathcal{N}_{\varepsilon}$, by Lemma \ref{lemma:C2-1-p-sum-C2-1}(i), we know $\Vert u_t-1\Vert_{C^2}<O(\varepsilon)$.
		If we denote
		\begin{align*}
			V(t)=\int_{\mathbb{S}^n}u_t\mathcal{Q}\left(\tau[u_t][k],I\right)dx,
		\end{align*}
		then $V(0)=V_k(K), V(1)=V_k(L)$. Using \eqref{equ:int-uQ-sysmetry}, we have
		\begin{align}
		V'(t)=&\int_{\mathbb{S}^n}v\mathcal{Q}\left(\tau[u_t][k],I\right)dx
			+k\int_{\mathbb{S}^n}u_t\mathcal{Q}\left(\tau[v][1],\tau[u_t][k-1],I\right)dx\nonumber
			\\
			\overset{\eqref{equ:int-uQ-sysmetry}}{=}&(k+1)\int_{\mathbb{S}^n}v\mathcal{Q}\left(\tau[u_t][k],I\right)dx,
			\label{equ:pf-Lemma 3.2}
		\end{align}
		and
		\begin{align*}
		V''(t)=&k(k+1)\int_{\mathbb{S}^n}v\mathcal{Q}\left(\tau[v][1],\tau[u_t][k-1],I\right)dx.
		\end{align*}
		From Taylor expansion (with integral form of the remainder) and $V_k(K)=V_k(L)$, it follows that
		\begin{align}
			0=&V(1)-V(0)=V'(0)+\int_0^1(1-t)V''(t)dt\nonumber
			\\
			=&(k+1)\int_{\mathbb{S}^n}v\mathcal{Q}\left(\tau[u_K][k],I\right)dx
			\nonumber
			\\&+\int_0^1\int_{\mathbb{S}^n}(1-t)k(k+1)v\mathcal{Q}\left(\tau[v][1],\tau[u_t][k-1],I\right)dxdt\nonumber
			\\
			=&(k+1)\int_{\mathbb{S}^n}vE_k(\lambda_K)dx
			+\int_0^1(1-t)k(k+1)B(t)dt,\label{equ:pf-int vEk<OE-1}
		\end{align}
		where we denote $B(t):=\int_{\mathbb{S}^n}v\mathcal{Q}\left(\tau[v][1],\tau[u_t][k-1],I\right)dx$.
		One can check that  $\tau[u_t]=\tau[u_t-1]+I$. Then we can write $B(t)$ as
		\begin{align}
			\label{equ:pf-int vEk<OE-2}
			B(t)=&\int_{\mathbb{S}^n}v\sum_{s=0}^{k-1}\binom{k-1}{s}\mathcal{Q}\left(\tau[v][1],\tau[u_t-1][s],I\right)dx.
		\end{align}
For index $s=0$, we denote $B_0$ by
	\begin{align}\label{equ:pf-int vEk<OE-3}
		B_0:=&\int_{\mathbb{S}^n}v\mathcal{Q}\left(\tau[v][1],I\right)dx\nonumber
		\nonumber\\
		=&\frac{1}{n}\int_{\mathbb{S}^n}v\left(\Delta v+nv\right)dx\nonumber
		\\=&\int_{\mathbb{S}^n}\left(-\frac{1}{n}|\nabla v|^2+v^2\right) dx=-\frac{1}{n}\Vert\nabla v\Vert_{L^2}^2+\Vert v\Vert_{L^2}^2,
	\end{align}
where we used integration by parts in the third equality. Moreover, we used integration by parts again to calculate $B(t)-B_0$ as follows
\begin{align}
	B(t)-B_0=&\int_{\mathbb{S}^n}v\sum_{s=1}^{k-1}\binom{k-1}{s}\mathcal{Q}\left(\tau[v][1],\tau[u_t-1][s],I\right)dx\nonumber
	\\
	\overset{\eqref{equ:Q=Qij*Aij}}{=}&\int_{\mathbb{S}^n} v\sum_{s=1}^{k-1}\binom{k-1}{s}\tau_{ij}[v]\frac{\partial\mathcal{Q}\left(\tau[v][1],\tau[u_t-1][s],I\right)}{\partial \tau_{ij}[v]}dx\nonumber
	\\
	\overset{\eqref{equ:Q_ij,j=0}}{=}&\int_{\mathbb{S}^n} \sum_{s=1}^{k-1}\binom{k-1}{s}\left(-v_jv_i+v^2\delta_{ij}\right)\frac{\partial\mathcal{Q}\left(\tau[v][1],\tau[u_t-1][s],I\right)}{\partial \tau_{ij}[v]}dx.\nonumber
\end{align}
Since $\Vert u_t-1\Vert_{C^2}<O(\varepsilon)$, we know $\left|\tau_{ij}[u_t-1]\right|\leq O(\varepsilon)$. When $s\geq 1$, we have
\begin{align*}
	\left|\frac{\partial\mathcal{Q}\left(\tau[v][1],\tau[u_t-1][s],I\right)}{\partial \tau_{ij}[v]}\right|\leq O(\varepsilon).
\end{align*}
Then we can estimate $|B(t)-B_0|$ as follows
\begin{align*}
	|B(t)-B_0|\leq O(\varepsilon)\Vert\nabla v\Vert_{L^2}^2+O(\varepsilon)\Vert v\Vert_{L^2}^2.
\end{align*}
Then
\begin{align}\label{equ:pf-int vEk<OE-4}
	B(t)-B_0= O(\varepsilon)\Vert\nabla v\Vert_{L^2}^2+O(\varepsilon)\Vert v\Vert_{L^2}^2.
\end{align}
Putting \eqref{equ:pf-int vEk<OE-3} and \eqref{equ:pf-int vEk<OE-4} into \eqref{equ:pf-int vEk<OE-1} yields
\begin{align*}
   \int_{\mathbb{S}^n}vE_k(\lambda_K)dx
	=&-k\int_0^1(1-t)B(t)dt
	\\
	=&-kB_0\int_0^1(1-t)dt-k\int_0^1(1-t)\left(B(t)-B_0\right)dt
	\\
	=&-\frac{k}{2}\left(-\frac{1}{n}\Vert\nabla v\Vert_{L^2}^2+\Vert v\Vert_{L^2}^2\right)+O(\varepsilon)\Vert\nabla v\Vert_{L^2}^2+O(\varepsilon)\Vert v\Vert_{L^2}^2.
\end{align*}
Then we complete the proof.
\end{proof}
	
	Consequently, we can express the deficit of $L^p$-Minkowski inequality for volume and intrinsic volumes in the following lemma.
	\begin{lemma}
		\label{lemma:I_k-w-barw-O(e)}
		Let $1\leq k\leq n$, and $K,L$ be    convex bodies  with support function $u_K(x)=1+w(x)\in C^3(\mathbb{S}^n)$ and  $u_L(x)=1+\bar{w}(x)\in C^3(\mathbb{S}^n)$, respectively. If $\Vert \bar{w} \Vert_{C^2}<\varepsilon$, $\Vert w \Vert_{C^2}<\varepsilon$ and $V_k(K)=V_k(L)$, then we have for $p=0$,
		\begin{align}
			\label{equ:deficit=p=0}
			\int_{\mathbb{S}^n}\left(\log(\frac{u_L}{u_K})\right)u_KE_k(\lambda_K)~dx
			= -\frac{k+1}{2}\Vert v\Vert_{L^2}^2
			+\frac{k}{2n}\Vert\nabla v\Vert_{L^2}^2+\mathscr{E},
		\end{align}
		and for $p\neq 0$,
		\begin{align}
			\label{equ:deficit=p!=0}
			\frac{1}{p}\left(\int_{\mathbb{S}^n}
			u_L^pu_K^{1-p}E_k(\lambda_K)~dx-V_k(K)\right)
			= -\frac{k+1-p}{2}\Vert v\Vert_{L^2}^2
			+\frac{k}{2n}\Vert\nabla v\Vert_{L^2}^2+\mathscr{E}.
		\end{align}
	\end{lemma}
	\begin{proof}
	  We take $\varepsilon$ small enough so that $\tau[u_K]>0,\tau[u_L]>0$, $\frac{1}{2}<u_K<\frac{3}{2}$, and $\left|\frac{v}{u_K}\right|<1$.
		By  Taylor expansion, we have
		 \begin{align*}
			\log\left(\frac{u_L}{u_K}\right)=&\log\left(1+\frac{v}{u_K}\right)
		\\=&\frac{v}{u_K}-\frac{1}{2}\left(\frac{v}{u_K}\right)^2+\sum_{i=3}^{\infty}(-1)^{i+1}\frac{1}{i}\left(\frac{v}{u_K}\right)^i
		\\=&\frac{v}{u_K}-\frac{1}{2}\left(\frac{v}{u_K}\right)^2+O(\varepsilon)\left(\frac{v}{u_K}\right)^2.
		\end{align*}
		The last equality above follows from the estimate $$\left|\sum_{i=3}^{\infty}(-1)^{i+1}\frac{1}{i}\left(\frac{v}{u_K}\right)^i\right|\leq O(\varepsilon)\left(\frac{v}{u_K}\right)^2.$$
		Together with Proposition \ref{lemma:int vEk<OE}, we obtain
		\begin{align*}       &\int_{\mathbb{S}^n}\left(\log(\frac{u_L}{u_K})\right)u_KE_k(\lambda_K)~dx
			\\
			=&\int_{\mathbb{S}^n}
			vE_k(\lambda_K)dx-\frac{1}{2}\int_{\mathbb{S}^n}
			v^2u_K^{-1}E_k(\lambda_K)dx+O(\varepsilon)\int_{\mathbb{S}^n}
			v^2u_K^{-1}E_k(\lambda_K)dx
			\\
			=&\frac{k}{2}\left(\frac{1}{n}\Vert\nabla v\Vert_{L^2}^2-\Vert v\Vert_{L^2}^2\right)+O(\varepsilon)\Vert\nabla v\Vert_{L^2}^2+O(\varepsilon)\Vert v\Vert_{L^2}^2
			-\frac{1}{2}\Vert v\Vert_{L^2}^2\left(1+O(\varepsilon)\right)+O(\varepsilon)\Vert v\Vert_{L^2}^2
			\\
			=&
			-\frac{k+1}{2}\Vert v\Vert_{L^2}^2
			+\frac{k}{2n}\Vert\nabla v\Vert_{L^2}^2+O(\varepsilon)\Vert  v\Vert_{L^2}^2+O(\varepsilon)\Vert \nabla v\Vert_{L^2}^2,
		\end{align*}
		where we use the fact that $u_K^{-1}=\frac{1}{1+O(\varepsilon)}=1+O(\varepsilon)$	and $E_{k}(\lambda_k)=1+O(\varepsilon)$ (see \eqref{equ:E_k(lambda)}).
		
		When $p\neq 0$, from  Taylor expansion, we have
		\begin{align*}
			\left(\frac{u_L}{u_K}\right)^p=&\left(1+\left(\frac{v}{u_K}\right)\right)^p
			\\
			=&1+p\left(\frac{v}{u_K}\right)+\frac{p(p-1)}{2}\left(\frac{v}{u_K}\right)^2
			+\sum_{i=3}^{\infty}\binom{p}{i}\left(\frac{v}{u_K}\right)^i
			\\
			=&1+p\left(\frac{v}{u_K}\right)+\frac{p(p-1)}{2}\left(\frac{v}{u_K}\right)^2+O(\varepsilon)\left(\frac{v}{u_K}\right)^2,
		\end{align*}
		where the last equality above follows from the estimate $$\left|\sum_{i=3}^{\infty}\binom{p}{i}\left(\frac{v}{u_K}\right)^i\right|\leq O(\varepsilon)\left(\frac{v}{u_K}\right)^2.$$
		Here $\binom{p}{i} := \frac{p(p-1)\cdots(p-i+1)}{i!}$ denotes the generalized binomial coefficient, for any real $p$ and integer $i \ge 0$.
		Together with Proposition \ref{lemma:int vEk<OE}, we obtain
		\begin{align*}
			&\frac{1}{p}\left(\int_{\mathbb{S}^n}
			u_L^pu_K^{1-p}E_k(\lambda_K)~dx-V_k(K)\right)
			\\
			=&\frac{1}{p}\int_{\mathbb{S}^n}
			u_K\left(\left(1+\left(\frac{v}{u_K}\right)\right)^p-1\right)E_k(\lambda_K)
			~dx
			\\
			=&\frac{1}{p}\int_{\mathbb{S}^n}
			u_K\left(p\left(\frac{v}{u_K}\right)+\frac{p(p-1)}{2}\left(\frac{v}{u_K}\right)^2+O(\varepsilon)\left(\frac{v}{u_K}\right)^2\right)E_k(\lambda_K)
			~dx
			\\
			=&
			-\frac{k+1-p}{2}\Vert v\Vert_{L^2}^2
			+\frac{k}{2n}\Vert\nabla v\Vert_{L^2}^2+O(\varepsilon)\Vert  v\Vert_{L^2}^2+O(\varepsilon)\Vert \nabla v\Vert_{L^2}^2.
		\end{align*}
		Then we complete the proof.
	\end{proof}
	
	Inspired by \cite{Cicalese-Leonardi-2012}*{eq. (73)}, we derive the lower bound on $\Vert v\Vert_{L^2}^2$.
	\begin{lemma}
		\label{lemma:DU>(2n+1)u}
		Let $p\in \mathbb{R}$,  $1\leq k\leq n$, and $K,L\in \mathcal{N}_{\varepsilon}$ be    convex bodies  with support function $u_K=1+w,u_L=1+\bar{w}\in C^3(\mathbb{S}^n)$, respectively. Assume $V_k(K)=V_k(L)$  and the $p$-weighted  Steiner points satisfy  $\mathcal{S}_p(K)=\mathcal{S}_p(L) $.
		Then we have
		\begin{align}
			\label{equ:equ:DU>(2n+1)u}
			\Vert\nabla(u_K-u_L)\Vert_{L^2}^2\geq2(n+1)\Vert u_K-u_L\Vert_{L^2}^2+\mathscr{E}.
		\end{align}
	\end{lemma}
	\begin{proof}
		Denote by $\{Y_{\ell}\}$ the spherical harmonics which form an orthonormal basis for $L^2(\mathbb{S}^n)$, i.e., 
		\begin{align}\label{equ:yiyj}
			\int_{\mathbb{S}^n}Y_{k}Y_{\ell} dx=\delta_{k \ell}, \quad \text{for\ } k,\ell=0,1,2,\cdots,
		\end{align}
		and
			\begin{align}
			\label{equ:Laplave Y=lambda Y}
			-\Delta_{\mathbb{S}^n} Y_{\ell}	=\lambda_{\ell}Y_{\ell}, \quad \text{for\ }\ell =0,1,2,\cdots,
		\end{align}
			where $\lambda_{\ell}=\ell(\ell+n-1)$ are the eigenvalues of the negative Laplace-Beltrami operator $-\Delta_{\mathbb{S}^n}$ on the sphere, with corresponding eigenfunctions being the spherical harmonics of degree $\ell$. 
		We write
		\begin{align}
			\label{equ:w-w-}
			w-\bar{w}=\sum_{\ell=0}^{\infty}a_{\ell}Y_{\ell}.
		\end{align}
		Recall that
		\begin{align*}
			Y_0=\frac{1}{\sqrt{|\mathbb{S}^n|}},\quad Y_1(x)=\<x,e\>,
		\end{align*}
		for a suitably chosen constant vector $e$. Here $|\mathbb{S}^n|$ denotes the surface area of $\mathbb{S}^n$.
		
		We have
		\begin{align*}
			a_0=\int_{\mathbb{S}^n}( w-\bar{w})Y_0~dx=\frac{1}{\sqrt{|\mathbb{S}^n|}}\int_{\mathbb{S}^n}( w-\bar{w})~dx,
		\end{align*}
		and
		\begin{align*}
			a_1=\int_{\mathbb{S}^n}( w-\bar{w})Y_1~dx=\left\langle \mathcal{S}(K)-\mathcal{S}(L),e\right\rangle. 
		\end{align*}
		Since $V_k(K)=V_k(L)$, from Lemma \ref{lemma:int vEk<OE} and $|E_k(\lambda_K)-1|\leq O(\varepsilon)$, we know
		\begin{align*}
			-{\sqrt{|\mathbb{S}^n|}}a_0=&\int_{\mathbb{S}^n}v dx=\int_{\mathbb{S}^n}vE_k(\lambda_K) dx+\int_{\mathbb{S}^n}v\left(1-E_k(\lambda_K)\right) dx
			\\
			=&-\frac{k}{2}\Vert v\Vert_{L^2}^2+\frac{k}{2n}\Vert \nabla v\Vert_{L^2}^2+\mathscr{E}+O(\varepsilon)\Vert v\Vert_{L^2}.
		\end{align*}
		This implies 
		\begin{align*}
			a_0=&\frac{1}{\sqrt{|\mathbb{S}^n|}}\left(\frac{k}{2}\Vert v\Vert_{L^2}^2-\frac{k}{2n}\Vert \nabla v\Vert_{L^2}^2+\mathscr{E}+O(\varepsilon)\Vert v\Vert_{L^2}\right)
			\\
			=&\frac{k}{2\sqrt{|\mathbb{S}^n|}}\Vert v\Vert_{L^2}^2-\frac{k}{2n\sqrt{|\mathbb{S}^n|}}\Vert \nabla v\Vert_{L^2}^2+\mathscr{E}+O(\varepsilon)\Vert v\Vert_{L^2}
			\\
			=&O(\varepsilon)\Vert v\Vert_{L^2}.
		\end{align*}
		Then
		\begin{align}\label{equ:a_0=OE}
			a_0^2=O(\varepsilon)\Vert v\Vert_{L^2}a_0=O(\varepsilon)\Vert v\Vert_{L^2}\left(\frac{k}{2\sqrt{|\mathbb{S}^n|}}\Vert v\Vert_{L^2}^2-\frac{k}{2n\sqrt{|\mathbb{S}^n|}}\Vert \nabla v\Vert_{L^2}^2+\mathscr{E}+O(\varepsilon)\Vert v\Vert_{L^2}\right)=\mathscr{E}.
		\end{align}
			For any $p\in\mathbb R$, set $\Phi_p(z)=(z^p-1)/p$ for $p\ne0$
		and $\Phi_0(z)=\log z$. Since $\Phi_p'(z)=z^{p-1}$ for $z>0$, we have
		\begin{equation*}
			\Phi_p(u_L)-\Phi_p(u_K)=vB,\qquad
			B(x)=\int_0^1 u_t(x)^{p-1}\,dt.
		\end{equation*}
		By the definition of the $p$-weighted Steiner point, $\mathcal{S}_p(K)=\mathcal{S}_p(L)$ gives $\int_{\mathbb S^n}(\Phi_p(u_L)-\Phi_p(u_K))x\,dx=0$.
		Thus, also in the logarithmic case,
		\begin{equation*}
			\int_{\mathbb S^n}vB Y_1\,dx=0.
		\end{equation*}
		As $u_t=1+O(\varepsilon)$ uniformly, $\|1-B\|_{L^\infty}\leq
		C_p\varepsilon$. Then we have
		\begin{equation*}
			|a_1|=\left| \int_{\mathbb{S}^n}vY_1dx\right| =\left| \int_{\mathbb{S}^n}v(1-B)Y_1dx\right|
			\leq C_p\varepsilon\|v\|_{L^2}.
		\end{equation*}
		For any $p\in\mathbb{R}$,  we have
		\begin{align}\label{equ:a_1<OE}
			a_1^2\leq O(\varepsilon)\Vert w-\bar{w}\Vert_{L^2}^2.
		\end{align}

		By \eqref{equ:yiyj} and  \eqref{equ:w-w-} we have  $$\Vert w-\bar{w}\Vert_{L^2}^2=\sum_{\ell=0}^{\infty}a_{\ell}^2.$$  
	From \eqref{equ:yiyj}$\thicksim$\eqref{equ:w-w-} and  integration by parts, we can calculate that
	\begin{align*}
		\Vert\nabla( w-\bar{w})\Vert_{L^2}^2=&-\int_{\mathbb{S}^n}(w-\bar{w})\cdot \Delta_{\mathbb{S}^n}(w-\bar{w}) dx
		\\
		\overset{\eqref{equ:w-w-}}{=}&-\int_{\mathbb{S}^n}
		\left(\sum_{\ell=0}^{\infty}a_{\ell}Y_{\ell}\right)\cdot
		\left( \sum_{k=0}^{\infty}a_{k}\Delta_{\mathbb{S}^n}Y_{k}\right) dx
		\\
		\overset{\eqref{equ:Laplave Y=lambda Y}}{=}&\int_{\mathbb{S}^n}
		\left(\sum_{\ell=0}^{\infty}a_{\ell}Y_{\ell}\right)\cdot
		\left( \sum_{k=0}^{\infty}a_k\lambda_{k}Y_{k}\right) dx
		\\
	\overset{\eqref{equ:yiyj}}{=}& \sum_{\ell=0}^{\infty}\lambda_{\ell}a_{\ell}^2.
	\end{align*}
		Noting that $\lambda_{\ell}=\ell(\ell+n-1)\geq 2(n+1)$ when $\ell\geq 2$, we derive
		\begin{align*}
			\Vert\nabla( w-\bar{w})\Vert_{L^2}^2=&\sum_{\ell=0}^{\infty}\lambda_{\ell}a_{\ell}^2    =\sum_{\ell=2}^{\infty}\lambda_{\ell}a_{\ell}^2+na_1^2
			\\
			\geq& 2(n+1)\sum_{\ell=2}^{\infty}a_{\ell}^2+na_1^2
			\\
			=&2(n+1)\sum_{\ell=0}^{\infty}a_{\ell}^2-2(n+1)a_0^2-(n+2)a_1^2
			\\
			\geq&2(n+1)\Vert w-\bar{w}\Vert_{L^2}^2+\mathscr{E}.
		\end{align*}
		In the last equality, we used \eqref{equ:a_0=OE} and \eqref{equ:a_1<OE}.
	\end{proof}

	\section{Stability of $L^p$-Minkowski inequality}\label{sec 4}
	Based on the propositions in Section \ref{sec 3}, this section aims to establish the stability of $L^p$-Minkowski inequality in terms of several distances. We shall prove Theorems \ref{thm:p-M-L2-Vk=Vk} and \ref{thm:p-M-index-kneqn}, 
	as well as Corollaries \ref{thm:p-M-L2} and \ref{thm:p-M-index-k!=n=!Vk}.
	\subsection{In terms of $L^2$ distance}
	We first prove the following theorem, which is equivalent to Corollary \ref{thm:p-M-L2}.
	\begin{theorem}
		\label{thm:p-M-L2-Vk=Vk}
		Let $K, L \in \mathcal{K}_o^{n+1}$ be $C^3$ convex bodies with $V_k(K)=V_k(L)$, $1\leq k\leq n$.  For all $\eta>0$, there exists $\varepsilon>0$  such that  if 
		$K,L\in\mathcal{N}_{\varepsilon}$ and  satisfy  $\mathcal{S}_p(K)=\mathcal{S}_p(L)$, then we have for $p=0$,
		\begin{align}\label{equ:p-M-Vk=p=0}
			\int_{\mathbb{S}^n}\left(\log(  \frac{u_L}{u_K})\right)u_KE_k(\lambda_K) ~dx
			\geq \left(\frac{1}{2}(k-1+\frac{2k}{n})-\eta\right)\delta_2^2(K,L);
		\end{align}
		and for $p\in (-k+1-\frac{2k}{n},0)\cup(0,+\infty)$,
		\begin{align}\label{equ:p-M-Vk=p!=0}
			\frac{1}{p}\left(\int_{\mathbb{S}^n}u_L^pu_K^{1-p}E_k(\lambda_K)~dx-V_k(K)\right)
			\geq \left(\frac{1}{2}(p+k-1+\frac{2k}{n})-{\eta}\right)\delta_2^2(K,L).
		\end{align}
	\end{theorem}
	
	\begin{proof}
		Substituting  \eqref{equ:equ:DU>(2n+1)u} into \eqref{equ:deficit=p=0} yields
		\begin{align}
			&\int_{\mathbb{S}^n}\left(\log(\frac{u_L}{u_K})\right)u_KE_k(\lambda_K)~dx		\nonumber
			\\
			=&\left(-\frac{k+1}{2}+O(\varepsilon)\right)\Vert\bar{w}-w\Vert_{L^2}^2
			+\frac{2(n+1)k}{2n}\left(\frac{1}{2(n+1)}+O(\varepsilon)\right)\Vert\nabla\bar{w}-\nabla w\Vert_{L^2}^2\nonumber
			\\
			\geq& \left(-\frac{k+1}{2}+O(\varepsilon)\right)\Vert\bar{w}-w\Vert_{L^2}^2
			+\frac{(n+1)k}{n}\left(1+O(\varepsilon)\right)\Vert\bar{w}- w\Vert_{L^2}^2\nonumber
			\\
			=&\left(\frac{1}{2}(k-1+\frac{2k}{n})+O(\varepsilon)\right)\delta_2^2(K,L),\nonumber
		\end{align}
	which proves \eqref{equ:p-M-Vk=p=0}.
	
	  For $p >-k+1-\frac{2k}{n}$ and $p\neq 0$, the analysis proceeds analogously to the case $p=0$. Putting \eqref{equ:equ:DU>(2n+1)u} into \eqref{equ:deficit=p!=0} yields
		\begin{align}
			&\frac{1}{p}\left(\int_{\mathbb{S}^n}u_L^pu_K^{1-p}E_k(\lambda_K)~dx-V_k(K)\right)\nonumber
			\\
			=&\left(-\frac{k+1-p}{2}+O(\varepsilon)\right)\Vert\bar{w}-w\Vert_{L^2}^2
			+\frac{2(n+1)k}{2n}\left(\frac{1}{2(n+1)}+O(\varepsilon)\right)\Vert\nabla\bar{w}-\nabla w\Vert_{L^2}^2\nonumber
			\\
			\geq&\left(-\frac{k+1-p}{2}+O(\varepsilon)\right)\Vert\bar{w}-w\Vert_{L^2}^2
			+\frac{(n+1)k}{n}\left(1+O(\varepsilon)\right)\Vert \bar{w}-  w\Vert_{L^2}^2\nonumber
			\\
			=& \left(\frac{1}{2}(p+k-1+\frac{2k}{n})+O(\varepsilon)\right)\delta_2^2(K,L),\nonumber
		\end{align}
		which proves \eqref{equ:p-M-Vk=p!=0}.
	\end{proof}
	Applying Theorem \ref{thm:p-M-L2-Vk=Vk} now yields Corollary \ref{thm:p-M-L2} for the unnormalized case (i.e., without requiring $V_k(K)=V_k(L)$).
	\begin{corollary}
		\label{thm:p-M-L2}
		Let $K, L \in \mathcal{K}_o^{n+1}$ be $C^3$ convex bodies and $1\leq k\leq n$. 
		For all $\eta>0$, there exists $\varepsilon>0$  such that if 
		$K,L\in\mathcal{N}_{\varepsilon}$, and   satisfy  $\mathcal{S}_p(K)=\mathcal{S}_p(L)=o$, then we have for $p=0$,
		\begin{align}\label{equ:p-M-L2-p=1}
			\int_{\mathbb{S}^n}\left(\log(\frac{u_L}{u_K})\right)u_KE_k(\lambda_K) ~dx
			-\frac{1}{k+1}V_k(K)\log\(\frac{V_k(L)}{V_k(K)}\)\nonumber
			\\
			\geq N^{k-1}\left(\frac{1}{2}(k-1+\frac{2k}{n})-\eta\right)\delta_2^2( K,NL);
		\end{align}
		and for $p\in (-k+1-\frac{2k}{n},0)\cup(0,+\infty)$, 
		\begin{align}\label{equ:p-M-L2-p=!1}
			\frac{1}{p}\left(\int_{\mathbb{S}^n}u_L^pu_K^{1-p}E_k(\lambda_K)~dx-V_k(K)^{\frac{k+1-p}{k+1}}V_k(L)^{\frac{p}{k+1}}\right)\nonumber
			\\
			\geq N^{k-1-p}\left(\frac{1}{2}(p+k-1+\frac{2k}{n})-{\eta}\right)\delta_2^2( K,NL),
		\end{align}
		where 
		$N=\left(\frac{V_k(K)}{V_k(L)}\right)^{\frac{1}{k+1}}$.
	\end{corollary}

	\begin{proof}
		For $K$ and $L$, there exists a constant $N>0$, such that $\widetilde{K}:=\frac{1}{N}K$ satisfies $V_k(\widetilde{K})=V_k(L)$. Then we have
		$$\left(\frac{1}{N}\right)^{k+1}V_k(K)=V_k(L),$$
		which implies
		\begin{align*}
			N=\left(\frac{V_k(K)}{V_k(L)}\right)^{\frac{1}{k+1}}.
		\end{align*}
		From \eqref{equ:Vk}, we have $V_k(K)=|\mathbb{S}^{n}|+O(\varepsilon)$, and $V_k(L)=|\mathbb{S}^{n}|+O(\varepsilon)$, then $N=1+O(\varepsilon)$. Denote by $u_{\widetilde{K}}=\frac{u_K}{N}:=1+\tilde{w}$  the support function of $\widetilde{K}$, then we have $\Vert \tilde{w}\Vert_{C^2}=O(\varepsilon)$.
		One can check that   $(\widetilde{K},L)$ satisfies $\mathcal{S}_p(\widetilde{K})=\mathcal{S}_p(L)=o$.

		By $K=N\widetilde{K}$ and Theorem \ref{thm:p-M-L2-Vk=Vk}, we have
		\begin{align*}
			&\int_{\mathbb{S}^n}\left(\log(\frac{u_L}{u_K})\right)u_KE_k(\lambda_K) ~dx
			-\frac{1}{k+1}V_k(K)\log\frac{V_k(L)}{V_k(K)}
			\\
			=&N^{k+1}\left(\int_{\mathbb{S}^n}\left(\log(\frac{u_L}{Nu_{\widetilde{K}}})\right)u_{\widetilde{K}}E_k(\lambda_{\widetilde{K}}) ~dx-\frac{1}{k+1}V_k(\widetilde{K})\log\frac{V_k(L)}{N^{k+1}V_k(\widetilde{K})}\right)
			\\
			=&N^{k+1}\int_{\mathbb{S}^n}\left(\log(\frac{u_L}{u_{\widetilde{K}}})\right)u_{\widetilde{K}}E_k(\lambda_{\widetilde{K}}) ~dx
			\\
			\geq& N^{k+1}\left(\frac{1}{2}(k-1+\frac{2k}{n})-\eta\right)\delta_2^2(\widetilde{K},L)
			.
		\end{align*}
		Since
		\begin{align*}
			\delta_2^2(\widetilde{K},L)=&\int_{\mathbb{S}^n}(u_{\widetilde{K}}-u_L)^2~dx
			=\int_{\mathbb{S}^n}(\frac{1}{N}u_K-u_L)^2~dx
			\\
			=&N^{-2}\int_{\mathbb{S}^n}(u_K-Nu_L)^2~dx
			\\
			=&N^{-2}\delta_2^2({K},NL),
		\end{align*}
we obtain the inequality \eqref{equ:p-M-L2-p=1}.
		
		Similarly, when $p> 0$ or $-k+1-\frac{2k}{n}<p<0$, 	by $K=N\widetilde{K}$ and Theorem \ref{thm:p-M-L2-Vk=Vk}, we have
		\begin{align*}
			&\frac{1}{p}\left(\int_{\mathbb{S}^n}u_L^pu_K^{1-p}E_k(\lambda_K)~dx-V_k(K)^{\frac{k+1-p}{k+1}}V_k(L)^{\frac{p}{k+1}}\right)
			\\
			=&\frac{1}{p}\left(\int_{\mathbb{S}^n}u_L^pu_{(N\widetilde{K})}^{1-p}E_k(\lambda_{(N\widetilde{K})})~dx-V_k({N\widetilde{K}})^{\frac{k+1-p}{k+1}}V_k(L)^{\frac{p}{k+1}}\right)
			\\
			=&N^{k+1-p}\frac{1}{p}\left(\int_{\mathbb{S}^n}u_L^pu_{\widetilde{K}}^{1-p}E_k(\lambda_{\widetilde{K}})~dx-V_k({\widetilde{K}})^{\frac{k+1-p}{k+1}}V_k(L)^{\frac{p}{k+1}}\right)
			\\
			\geq& N^{k-1-p}\left(\frac{1}{2}(p+k-1+\frac{2k}{n})-{\eta}\right)\delta_2^2( K,NL).
		\end{align*}
	\end{proof}

	\subsection{In terms of relative asymmetry index}
	We can derive Theorem \ref{thm:p-M-index-kneqn} as follows, which constitute some equivalent formulations of Corollary \ref{thm:p-M-index-k!=n=!Vk}.
	\begin{theorem}
		\label{thm:p-M-index-kneqn}
		Let $K, L \in \mathcal{K}_o^{n+1}$ be $C^3$ convex bodies with $V_k(K)=V_k(L)$, $1\leq k\leq n$.
		For all $\eta>0$, there exists $\varepsilon>0$  such that if 
		$K,L\in\mathcal{N}_{\varepsilon}$ and $K,L$ satisfy  $\mathcal{S}_p(K)=\mathcal{S}_p(L)$, then we have for $p=0$,
		\begin{align}\label{equ:p-M-index-k!=n-p=0}
			\frac{\int_{\mathbb{S}^n}\log(\frac{u_L}{u_K})u_KE_k(\lambda_K)~dx}{V_{k}(K)}\geq \left(\frac{(k-1)n+2k}{2n(n+1)^2}-\eta\right)\alpha(K,L)^2;
		\end{align}	 	
		and for $p\in(1-k-\frac{2k}{n},0)\cup(0,+\infty)$,
		\begin{align}\label{equ:p-M-index-k!=n-p!=0}
			\frac{\int_{\mathbb{S}^n}u_L^pu_K^{1-p}E_{k}(\lambda_K)~dx-V_k(K)}{pV_k(K)}\geq \left( \frac{n(k-1+p)+2k}{2n(n+1)^2}-
			\eta\right)\alpha(K,L)^2.
		\end{align}
	\end{theorem}

	As a consequence of the above theorem, we deduce corollary that hold without the normalization condition $V_k(K)=V_k(L)$.
	\begin{corollary}
		\label{thm:p-M-index-k!=n=!Vk}
		Let $K, L \in \mathcal{K}_o^{n+1}$ be $C^3$ convex bodies and $1\leq  k\leq  n$.  For all $\eta>0$, there exists $\varepsilon>0$  such that if $K,L\in\mathcal{N}_{\varepsilon}$ and  satisfy  $\mathcal{S}_p(K)=\mathcal{S}_p(L)=o$, then we have for $p=0$,
		\begin{align}\label{equ:p-M-index-k!=n-p=0!=Vk}
			\int_{\mathbb{S}^n}\left(\log(\frac{u_L}{u_K})\right)u_KE_k(\lambda_K) ~dx
			-\frac{1}{k+1}V_k(K)\log\frac{V_k(L)}{V_k(K)}\nonumber
			\\
			\geq V_k(K) \left(\frac{(k-1)n+2k}{2n(n+1)^2}-\eta\right)\alpha(K,L)^2;
		\end{align}
	 and for $p\in (-k+1-\frac{2k}{n},0)\cup(0,+\infty)$, 
		\begin{align}\label{equ:p-M-index-k!=n-p!=0!=Vk}
			\frac{1}{p}\left(\int_{\mathbb{S}^n}u_L^pu_K^{1-p}E_k(\lambda_K)~dx-V_k(K)^{\frac{k+1-p}{k+1}}V_k(L)^{\frac{p}{k+1}}\right)\nonumber
			\\
			\geq N^{-p}V_k(K) \left(\frac{n(k-1+p)+2k}{2n(n+1)^2}-{\eta}\right)\alpha(K,L)^2.
		\end{align}
	\end{corollary}
	Before proving Theorem \ref{thm:p-M-index-kneqn}, we  study  the relative asymmetry index $\alpha(K,L)$  (defined by \eqref{equ:def-asymmetry}).

	\begin{lemma}
		Suppose $K, L \in \mathcal{K}_o^{n+1}$ are $C^2$ convex bodies, with support function $u_K$ and $u_L$, respectively.  Then we have
			\begin{align*}
			\left|K\Delta L\right|
			\leq \int_0^1\int_{\mathbb{S}^n}
			\left|u_K-u_L
			\right|
			\cdot \det\left(\nabla^2u_t+u_t I\right)
			dx dt
			.
		\end{align*}
		Here $u_t$ is the support function of $(1-t)K+tL$ ($t\in[0,1]$).
		Moreover, if $K, L \in\mathcal{N}_{\varepsilon}$, we have
		\begin{align}
			\left|K\Delta L\right|
			\leq\left(1+O(\varepsilon)\right)\cdot |\mathbb{S}^n|^\frac12\cdot \Vert u_K-u_L\Vert_{L^2}
			.\label{equ:pf-lemma-alpha>-4}
		\end{align}
	\end{lemma}
	\begin{proof}
			Let $\mathscr{A}_K^{*} : \mathbb{S}^n \to \mathbb{S}^n$ be the reverse radial Gauss map of convex body $K$ (for details, see \cite[eq. (1.5)]{Chen-Li-2021}). Let $\rho_K$ be the radial function of hypersurface $\partial K$, then we have $$\rho_K^2(\xi)=u_K^2(x)+|\nabla u_K|^2(x),$$ where $\xi=\mathscr{A}_K^{*}(x)$ or  $x(\xi)=\nu\left(\rho(\xi)\cdot \xi\right)\in\mathbb{S}^n$ for $x,\xi\in \mathbb{S}^n$ (see \cite{Chen-Li-2021}*{eq. (2.4)} or \cite{Ding-Li-2025-JDE}*{Section 2.3}),
		and the determinant of the Jacobian of $\mathscr{A}_K^{*}$ is given by \cite[eq. (2.28)]{Ding-Li-Trans}
		\begin{align}
			\label{equ:pf-lemma-alpha>-1}
			|\operatorname{Jac}\mathscr{A}_K^{*}|(x)=\left| \frac{d\xi}{dx}\right| =\frac{u_K(x)\sigma_n(\lambda_K(x))}{\rho_K^{n+1}\left(\mathscr{A}^{*}(x)\right)}.
		\end{align}
		Here $d\xi:=d\mu_{\mathbb{S}^n}(\xi)$ denotes the standard induced volume element on $\mathbb{S}^n$.
		We assume  $K_t=(1-t)K+tL$ ($t\in[0,1]$),  the support function and radial function  of $K_t$ are $u_t(x)$ and $\rho_t(\xi)$ respectively, with  $\xi=\mathscr{A}_{K_t}^{*}(x)$ and
		$u_t=(1-t)u_K+tu_L.$ By \cite{Ding-Li-Trans}*{eq. (2.26)}, we have
		\begin{align}
			\label{equ:pf-lemma-alpha>-2}
			\frac{\partial \rho_t(\xi)}{\partial t}=\frac{\rho_t}{u_t}\frac{\partial u_t(x)}{\partial t}=\frac{\rho_t}{u_t}\cdot v(x).
		\end{align}

		Note that
		\begin{align}
			\label{equ:|M-Delta-W|}
			\left|K\Delta L\right|=&\int_{\mathbb{S}^n} \int_{\min \{\rho_0,\rho_1\}}^{\max \{\rho_0,\rho_1\}}r^ndrd\xi=\frac{1}{n+1}\int_{\mathbb{S}^n}
			\left|
			\rho_1^{n+1}-{\rho}_0^{n+1}
			\right|
			d\xi,
		\end{align}
		where we use that,   $\max\{a^n,b^n\}=\left(\max \{a,b\}\right)^n$,
		$\min\{a^n,b^n\}=\left(\min \{a,b\}\right)^n$,
		$\max \{a,b\}=\frac{a+b+|a-b|}{2}$,
		and $\min \{a,b\}=\frac{a+b-|a-b|}{2}$
		for any $n\geq1$ and $a,b\in\mathbb{R}^+$.
		Putting \eqref{equ:pf-lemma-alpha>-1} and \eqref{equ:pf-lemma-alpha>-2} into \eqref{equ:|M-Delta-W|}, we have
		\begin{align}
			\left|K\Delta L\right|
			\leq &\int_0^1\int_{\mathbb{S}^n}
			\rho_t^n\cdot 	\left|\frac{\partial \rho_t}{\partial t}
			\right|
			d\xi dt\nonumber
			\\
			= &\int_0^1\int_{\mathbb{S}^n}
			\left|v(x)
			\right|
			\cdot \det\left(\nabla^2u_t+u_t I\right)
			dx dt. \label{equ:pf-lemma-alpha>-10}
		\end{align}
	The first assertion follows. If $K,L\in \mathcal{N}_{\varepsilon}$, we have $K_t:=(1-t)K+tL\in\mathcal{N}_{\varepsilon}$ with $t\in[0,1]$ (from Lemma \ref{lemma:C2-1-p-sum-C2-1} (i)), then we have $\tau_{ij}[u_t]=\delta_{ij}+O(\varepsilon)$ and
		\begin{align*}
			E_n\left(\nabla^2u_t+u_t I\right)=1+O(\varepsilon).
		\end{align*}
		Putting it into \eqref{equ:pf-lemma-alpha>-10} yields
		\begin{align*}
			\left|K\Delta L\right|\leq \left(1+O(\varepsilon)\right)\cdot |\mathbb{S}^n|^\frac12\cdot \Vert v\Vert_{L^2},
		\end{align*}
		where we used the H\"older inequality. This proves \eqref{equ:pf-lemma-alpha>-4}.
	\end{proof}
		
	\begin{lemma}
		\label{lemma:asymmetry}
		Suppose $K, L \in \mathcal{K}_o^{n+1}$ are $C^3$ convex bodies, with support function $u_K=1+w$ and $u_L=1+\bar{w}$, respectively.  If $\|w\|_{C^2}<\varepsilon$, $\|\bar{w}\|_{C^2}<\varepsilon$, and $V_k(K)=V_{k}(L)$, then we have
		\begin{align}
			\label{equ:lemma-asymmetry}
			\alpha(K,L)^2\leq \frac{(n+1)^2}{|\mathbb{S}^n|}\cdot
			\|w-\bar{w}\|_{L^2}^2
			+O(\varepsilon)\cdot
			\|w-\bar{w}\|_{L^2}^2.
		\end{align}
		\end{lemma}
	\begin{proof}
	
	Since $K,L\in \mathcal{N}_{\varepsilon}$, we have $K_t:=(1-t)K+tL\in\mathcal{N}_{\varepsilon}$ with $t\in[0,1]$ (from Lemma \ref{lemma:C2-1-p-sum-C2-1} (i)). Let $u_t$ be the support function of $K_t$,  we have $\tau_{ij}[u_t]=\delta_{ij}+O(\varepsilon)$ and
	\begin{align}
		\label{equ:pf-lemma-alpha>-5}
		E_k\left(\nabla^2u_t+u_t I\right)=1+O(\varepsilon).
	\end{align}
	Since $V_k(K)=V_k(L)$, from \eqref{equ:pf-Lemma 3.2}, we have
	\begin{align}
		\label{equ:pf-lemma-alpha>-6}
		0=V_k(L)-V_{k}(K)=(k+1)\int_0^1\int_{\mathbb{S}^n}vE_{k}(\tau[u_t])dxdt.
	\end{align}
	Then, we have
	\begin{align}
		\label{equ:pf-lemma-alpha>-7}
		\big||L|-|K|\big|=&\left|\int_0^1\int_{\mathbb{S}^n}vE_{n}(\tau[u_t])dxdt\right|\nonumber
		\\
		\overset{\eqref{equ:pf-lemma-alpha>-6}}{=}&\left|\int_0^1\int_{\mathbb{S}^n}v\left(E_{n}(\tau[u_t])-E_{k}(\tau[u_t])\right)dxdt\right|\nonumber
		\\
		\overset{\eqref{equ:pf-lemma-alpha>-5}}{\leq}&O(\varepsilon)\Vert v\Vert_{L^2},
	\end{align}
where we used the H\"older inequality again.
	Let $\W_{K}=rL$, where $r>0$ is a constant chosen such that $|K|=|\W_{K}|=r^{n+1}|L|$, then
	 \begin{align}\label{equ:pf-lemma-alpha>-8}
		r^{n+1}=\frac{|K|}{|L|}.
	\end{align}	
  	Combining  with  $|L|=\frac{|\mathbb{S}^n|}{(n+1)}+O(\varepsilon)$, and $|K|=\frac{|\mathbb{S}^n|}{(n+1)}+O(\varepsilon)$ we have $r=1+O(\varepsilon)$ and $rL\in\mathcal{N}_{\varepsilon}$. From \eqref{equ:pf-lemma-alpha>-8}, we have
  	\begin{align*}
  		r-1=\frac{|K|-|L|}{|L|\sum_{j=0}^nr^j},
  	\end{align*}
  	which implies
  	\begin{align}
  		\label{equ:pf-lemma-alpha>-9}
  		\left|
  		r-1\right|\leq O(\varepsilon)\Vert v\Vert_{L^2}.
  	\end{align}
   Then we have
   \begin{align*}
   	\Vert ru_L-u_K\Vert_{L^2}=&\Vert (r-1)u_L+v\Vert_{L^2}
   	\\
   	\leq &|r-1|\cdot \Vert u_L\Vert_{L^2}+\Vert v\Vert_{L^2}
   	\\
   	\leq& \left(1+O(\varepsilon)\right)\Vert v\Vert_{L^2}.
   \end{align*}
	Putting it into \eqref{equ:pf-lemma-alpha>-4} yields
	\begin{align*}
	|K\Delta\W_{K}|=|K\Delta(rL)|\leq \left(1+O(\varepsilon)\right)\cdot |\mathbb{S}^n|^\frac12\cdot \Vert v\Vert_{L^2}.
	\end{align*}
Since $|\W_K|=\frac{V_n(K)}{n+1}=\frac{|\mathbb{S}^n|}{n+1}(1+O(\varepsilon))$, 
we have
\begin{align*}
	\frac{|K\Delta\W_{K}|^2}{|\W_{K}|^2}
	\leq \left(\frac{(n+1)^2}{|\mathbb{S}^n|}
	+O(\varepsilon)\right)
	\|w-\bar{w}\|_{L^2}^2.
\end{align*}
On the other hand,
	\begin{align*}
		\alpha(K,L)=\min _{y \in \mathbb{R}^{n+1}}\left\{\frac{|K \Delta(r L+y)|}{|K|}:|r L|=|K|\right\}\leq \frac{|K \Delta \W_{K}|}{|\W_{K}|}.
	\end{align*}
 	Thus we obtain
		\begin{align*}
			\alpha(K,L)^2\leq \left(\frac{(n+1)^2}{|\mathbb{S}^n|}
			+O(\varepsilon)\right)
			\|w-\bar{w}\|_{L^2}^2.
		\end{align*}
	\end{proof}

	
	Now we can prove Theorem \ref{thm:p-M-index-kneqn} as follows.

	
	\begin{proof}[\textbf{Proof of Theorem \ref{thm:p-M-index-kneqn}}]
		Let  $1 \le k \le n$.	Combining Lemma \ref{lemma:I_k-w-barw-O(e)} with Lemma \ref{lemma:DU>(2n+1)u}, we obtain for $p=0$:
		\begin{align}   \label{equ:I_k-IkO>u-u-dudu}
			\int_{\mathbb{S}^n}\left(\log(\frac{u_L}{u_K})\right)u_KE_k(\lambda_K)~dx
			=&-\frac{k+1}{2}\Vert\bar{w}-w\Vert_{L^2}^2
			+\frac{k}{2n}\Vert\nabla\bar{w}-\nabla w\Vert_{L^2}^2+\mathscr{E}\nonumber
			\\
			\geq&(1+O(\varepsilon))\left(\frac{k}{2n}-\frac{k+1}{2}\frac{1}{2(n+1)}\right)\Vert\nabla\bar{w}-\nabla w\Vert_{L^2}^2\nonumber
			\\
			=&(1+O(\varepsilon))\frac{n(k-1)+2k}{4n(n+1)}\Vert\nabla\bar{w}-\nabla w\Vert_{L^2}^2.
		\end{align}
				Since $V_k(K)=V_k(L)$, from Lemma \ref{lemma:asymmetry} and Lemma \ref{lemma:DU>(2n+1)u}, we get
		\begin{align*}
			\alpha(K,L)^2\leq &\frac{(n+1)^2}{|\mathbb{S}^n|}\cdot
			(1+O(\varepsilon))\cdot \|w-\bar{w}\|_{L^2}^2
			\nonumber
			\\
			\leq &
			\frac{(n+1)^2}{|\mathbb{S}^n|}\frac{1}{2(n+1)}(1+O(\varepsilon))
			\|\nabla w-\nabla \bar{w}\|_{L^2}^2
			\nonumber
			\\
			=&\frac{n+1}{2|\mathbb{S}^n|}\cdot (1+O(\varepsilon))\cdot
			\|\nabla w-\nabla\bar{w}\|_{L^2}^2.
		\end{align*}
		Putting it into \eqref{equ:I_k-IkO>u-u-dudu}, we have
		\begin{align*}
			\int_{\mathbb{S}^n}\left(\log(\frac{u_L}{u_K})\right)u_KE_k(\lambda_K)~dx
			\geq &
			\frac{|\mathbb{S}^n|}{n+1}\cdot
			\frac{n(k-1)+2k}{2n(n+1)}\left(1+O(\varepsilon)\right)\cdot\alpha(K,L)^2\nonumber
			\\
			=&
			|\mathbb{S}^n|\cdot
			\frac{n(k-1)+2k}{2n(n+1)^2}\left(1+O(\varepsilon)\right)\cdot\alpha(K,L)^2			.
		\end{align*}
		By  \eqref{equ:Vk} we have  $V_k(K)=|\mathbb{S}^n|+O(\varepsilon)=(1+O(\varepsilon))|\mathbb{S}^n|$. This completes the proof of  \eqref{equ:p-M-index-k!=n-p=0}.
		
		For $p \in \left(1-k-\frac{2k}{n},0\right)\cup (0, \infty)$,
		by Lemma \ref{lemma:I_k-w-barw-O(e)} and Lemma \ref{lemma:DU>(2n+1)u}, we have
		\begin{align}\label{equ:proof-thm-index==1}
			&\frac{\int_{\mathbb{S}^n}u_L^pu_K^{1-p}E_{k}(\lambda_K)~dx-V_k(K)}{p}\nonumber
			\\
			=&-\frac{k+1-p}{2}\Vert\bar{w}-w\Vert_{L^2}^2
			+\frac{k}{2n}\Vert\nabla\bar{w}-\nabla w\Vert_{L^2}^2+\mathscr{E}\nonumber
			\\
			\geq& -\frac{k+1-p}{2}\Vert\bar{w}-w\Vert_{L^2}^2
			+\frac{k\cdot 2(n+1)}{2n}\Vert \bar{w}-  w\Vert_{L^2}^2+O(\varepsilon)\cdot \| w-\bar{w}\|_{L^2}^2\nonumber
			\\
			=&\frac{1}{2}\left(p-1+k+\frac{2k}{n}\right)\Vert\bar{w}-w\Vert_{L^2}^2
			+O(\varepsilon)\cdot \| w-\bar{w}\|_{L^2}^2.
		\end{align}
		Substituting \eqref{equ:lemma-asymmetry} into \eqref{equ:proof-thm-index==1} yields
		\begin{align*}
			\frac{\int_{\mathbb{S}^n}u_L^pu_K^{1-p}E_{k}(\lambda_K)~dx-V_k(K)}{p}\geq \frac{|\mathbb{S}^n|}{2(n+1)^2}\left(p-1+k+\frac{2k}{n}+O(\varepsilon)\right)
			\cdot\alpha(K,L)^2.
		\end{align*}
			Since $V_{k}(K)=(1+O(\varepsilon))|\mathbb{S}^n|$, we obtain \eqref{equ:p-M-index-k!=n-p!=0}.
	\end{proof}
	\begin{proof}[\textbf{Proof of Corollary \ref{thm:p-M-index-k!=n=!Vk}}] 	
		For $K$ and $L$, there exists a constant $N>0$, such that $\widetilde{K}:=\frac{1}{N}K$ satisfies $V_k(\widetilde{K})=V_k(L)$. Then we have
		$$\left(\frac{1}{N}\right)^{k+1}V_k(K)=V_k(L),$$
		which implies
		\begin{align*}
			N=\left(\frac{V_k(K)}{V_k(L)}\right)^{\frac{1}{k+1}}.
		\end{align*}
		From \eqref{equ:Vk}, we have $V_k(K)= |\mathbb{S}^{n}|+O(\varepsilon)$, and $V_k(L)=|\mathbb{S}^{n}|+O(\varepsilon)$, then $N=1+O(\varepsilon)$. Denote by $u_{\widetilde{K}}=\frac{u_K}{N}:=1+\tilde{w}$  the support function of $\widetilde{K}$, then we have $\Vert \tilde{w}\Vert_{C^2}=O(\varepsilon)$.  One can check that  pair $(\widetilde{K},L)$ satisfies $\mathcal{S}_p(\widetilde{K})=\mathcal{S}_p(L)$=o.
		
		By $K=N\widetilde{K}$ and Theorem \ref{thm:p-M-index-kneqn}, we have
		\begin{align*}
			&\int_{\mathbb{S}^n}\left(\log(\frac{u_L}{u_K})\right)u_KE_k(\lambda_K) ~dx
			-\frac{1}{k+1}V_k(K)\log\frac{V_k(L)}{V_k(K)}
			\\
			=&\int_{\mathbb{S}^n}\left(\log(\frac{u_L}{Nu_{\widetilde{K}}})\right)Nu_{\widetilde{K}}N^kE_k(\lambda_{\widetilde{K}}) ~dx-\frac{1}{k+1}N^{k+1}V_k(\widetilde{K})\log\frac{V_k(L)}{N^{k+1}V_k(\widetilde{K})}
			\\
			\geq& N^{k+1}\left(\frac{(k-1)n+2k}{2n(n+1)^2}-\eta\right)\alpha(L,\widetilde{K})^2V_k(\widetilde{K})
			\\
			=& {V_k(K)}\left(\frac{(k-1)n+2k}{2n(n+1)^2}-\eta\right)\alpha(L,\widetilde{K})^2
			.
		\end{align*}
		Since
		\begin{align*}
			\alpha(L,\widetilde{K})=&\min _{y \in \mathbb{R}^{n+1}}\left\{\frac{|\widetilde{K} \Delta(r L+y)|}{|\widetilde{K}|}:|r L|=|\widetilde{K}|\right\}
			\\
			=&\min _{y \in \mathbb{R}^{n+1}}\left\{\frac{|(\frac{1}{N}K) \Delta(\frac{1}{N}rN L+y)|}{|{\frac{1}{N}K}|}:|r L|=|\frac{1}{N}K|\right\}
			\\
			=&\min _{y \in \mathbb{R}^{n+1}}\left\{\frac{|K \Delta(rN L+y)|}{|{K}|}:|rN L|=|K|\right\}
			\\
			=&\alpha(K,L),
		\end{align*}
		we obtain the inequality \eqref{equ:p-M-index-k!=n-p=0!=Vk}.
		
		Similarly, we can derive \eqref{equ:p-M-index-k!=n-p!=0!=Vk} from Theorem \ref{thm:p-M-index-kneqn}.
	\end{proof}
	
	\section{Stability of $L^p$-Brunn--Minkowski inequality}\label{sec 5}
	In this section, we prove the stability of the  $L^p$-Brunn--Minkowski inequality stated in Theorem \ref{thm:p-BM-L2-Vk=Vk}. Moreover, we also derive the unnormalized versions (i.e., without requiring $V_k(K)=V_k(L)$) given in Theorems \ref{thm:p-BM-L2-Vk!=Vk} and \ref{thm:p-BM-L2-Vk!=Vk-index}. Before presenting the proofs, we first establish the following necessary lemma, which is inspired by \cite{Patsalos-Saroglou-arxiv}*{Lemma 4.1}.
	
	\begin{lemma}\label{lemma:C2-1-p-sum-C2-1}
		Assume $p\in\mathbb{R}$.	Given $\eta>0$, there exists a positive constant $\varepsilon$ such that, for every $K,L\in\mathcal{N}_{\varepsilon}$ and $t\in[0,1]$, the following hold.
		\begin{itemize}
			\item [(i)] $(1-t) \cdot K+_p t \cdot L \in \mathcal{N}_{\eta}$;
			\item [(ii)] $u_{(1-t) \cdot K+_p t \cdot L}=\left((1-t) u_K^p+t u_L^p\right)^{1 / p}$, \quad for $p\neq 0$;
			\item [(iii)] $u_{(1-t) \cdot K+_0 t \cdot L}= u_K^{1-t} u_L^{t}$; 
			\item [(iv)] $\left(\frac{V_j(\mathbb{B}^{n+1})}{V_j(K)}\right)^{\frac{1}{j+1}} K \in \mathcal{N}_{\eta}$.
		\end{itemize}
	\end{lemma}
	\begin{proof}
	We first consider the case $p\neq 0$.	Assume  $u_K=1+w$ and $u_L=1+\bar{w}$. Since $t$ belongs to the compact set $[0, 1]$, for  sufficiently small $\varepsilon$ and $K,L\in\mathcal{N}_{\varepsilon}$, we have
		\begin{align}
			\label{equ:pf-lemC2-1-p-sum-C2-1-1}
			\left(\left(1-t\right) u_{K}^p+t u_{L}^p\right)^{1 / p}=\left(\left(1-t\right) (1+O(\varepsilon))+t (1+O(\varepsilon))\right)^{1 / p}=1+O(\varepsilon),
		\end{align}
		and
		\begin{align}
			\label{equ:pf-lemC2-1-p-sum-C2-1-2}
			&\nabla_{e_i}\left(\left(1-t\right) u_{K}^p+t u_{L}^p\right)^{\frac{1}{p}} \nonumber		
			\\
			=&\frac{1}{p}\left(\left(1-t\right) u_{K}^p+t u_{L}^p\right)^{\frac{1-p}{p}}\left(\left(1-t\right) pu_{K}^{p-1}\nabla_{e_i}w+t pu_{L}^{p-1}\nabla_{e_i}\bar{w}\right)
			\nonumber\\
			=&O(\varepsilon).
		\end{align}
		Furthermore, we have
		\begin{align}
			\label{equ:pf-lemC2-1-p-sum-C2-1-3}
			&\nabla_{e_j}\nabla_{e_i}\left(\left(1-t\right) u_{K}^p+t u_{L}^p\right)^{\frac{1}{p}} \nonumber		
			\\
			=&\frac{1-p}{p^2}\left(\left(1-t\right) u_{K}^p+t u_{L}^p\right)^{\frac{1-2p}{p}}\left(\left(1-t\right) pu_{K}^{p-1}\nabla_{e_i}w+t pu_{L}^{p-1}\nabla_{e_i}\bar{w}\right)\cdot
			\nonumber\\
			&\left(\left(1-t\right) pu_{K}^{p-1}\nabla_{e_j}w+t pu_{L}^{p-1}\nabla_{e_j}\bar{w}\right)
			+
			\frac{1}{p}\left(\left(1-t\right) u_{K}^p+t u_{L}^p\right)^{\frac{1-p}{p}}\cdot\nonumber
			\\
			&\left(\left(1-t\right) p(p-1)u_{K}^{p-2}\nabla_{e_j}w\nabla_{e_i}w+t p(p-1)u_{L}^{p-2}\nabla_{e_j}\bar{w}\nabla_{e_i}\bar{w}\right)\nonumber
			\\
			&+
			\frac{1}{p}\left(\left(1-t\right) u_{K}^p+t u_{L}^p\right)^{\frac{1-p}{p}}\left(\left(1-t\right) pu_{K}^{p-1}\nabla_{e_j}\nabla_{e_i}w+t pu_{L}^{p-1}\nabla_{e_j}\nabla_{e_i}\bar{w}\right)\nonumber
			\\
			=&O(\varepsilon).
		\end{align}
			By \eqref{equ:pf-lemC2-1-p-sum-C2-1-1}$\thicksim$\eqref{equ:pf-lemC2-1-p-sum-C2-1-3}, we see $$\left\Vert \left(\left(1-t\right) u_{K}^p+t u_{L}^p\right)^{\frac{1}{p}} -1\right\Vert_{C^2}=O(\varepsilon).$$
		From \eqref{equ:pf-lemC2-1-p-sum-C2-1-1} and \eqref{equ:pf-lemC2-1-p-sum-C2-1-3}, it follows that the matrix $\left\{\tau_{ij}\left[\left(\left(1-t\right) u_{K}^p+t u_{L}^p\right)^{1 / p}\right]\right\}$ is positive definite for any $t\in[0,1]$ and $p\neq 0$. Then $(1+t )\cdot K+_p t\cdot L$ are convex bodies with support functions $u_{(1-t) \cdot K+_p t \cdot L}=\left((1-t) u_K^p+t u_L^p\right)^{1 / p}$, and $(1-t) \cdot K+_p t \cdot L\in\mathcal{N}_{O(\varepsilon)}$.
		This completes the proof of parts (i) for $p\neq 0$ and (ii).
		
		Similarly, for the case $p=0$, we can prove that
		$\left\Vert u_K^{1-t} u_L^{t}-1\right\Vert_{C^2}=O(\varepsilon)$ and matrix $\{ \tau_{ij}\left[u_K^{1-t} u_L^{t}\right] \}$ is positive definite for sufficiently small $\varepsilon$. This completes the proof of parts (i) for $p=0$ and (iii).
		
		From \eqref{equ:Vk}, we can check that $\left(\frac{V_j(\mathbb{B}^{n+1})}{V_j(K)}\right)^{1 / j}=1+O(\varepsilon)$, which completes the proof for  (iv).  	
	\end{proof}
	
	\begin{proof}[\textbf{Proof of Theorem \ref{thm:p-BM-L2-Vk=Vk}}]
		  We first consider the case $p=0$. From Lemma \ref{lemma:C2-1-p-sum-C2-1},  for  sufficiently small $\varepsilon$,  $K,L\in\mathcal{N}_{\varepsilon}$, and each $t\in[0,1]$, we have $Q_t:=(1-t)\cdot K+_0t\cdot L\in\mathcal{N}_{O(\varepsilon)}$ is a convex body with  support function $u_{Q_t}=u_K^{1-t} u_L^{t}$. Since $K$ and $L$ satisfy $\mathcal{S}_0(K)=\mathcal{S}_0(L)=o$, we know  $Q_t$ also satisfies $\mathcal{S}_0(Q_t)=o$.
		Then  from \eqref{equ:p-M-L2-p=1} we have
		\begin{align}
			\label{equ:pf-thm:p-BM-L2-Vk=Vk-4}
			&-t\int_{\mathbb{S}^n}\left(\log(\frac{u_L}{u_K})\right)u_{Q_t}E_k(\lambda_{Q_t}) ~dx=\int_{\mathbb{S}^n}\left(\log(\frac{u_K}{u_{Q_t}})\right)u_{{Q_t}}E_k(\lambda_{Q_t}) ~dx
			\nonumber
			\\
			\geq&\frac{1}{k+1}V_k(Q_t)\log\frac{V_k(K)}{V_k(Q_t)}
			+ N_1^{k-1}\left(\frac{1}{2}(k-1+\frac{2k}{n})-\eta\right)\delta_2^2(N_1 K,Q_t),
		\end{align}
		and
		\begin{align}
			\label{equ:pf-thm:p-BM-L2-Vk=Vk-5}
			&(1-t)\int_{\mathbb{S}^n}\left(\log(\frac{u_L}{u_K})\right)u_{Q_t}E_k(\lambda_{Q_t}) ~dx=\int_{\mathbb{S}^n}\left(\log(\frac{u_L}{u_{Q_t}})\right)u_{{Q_t}}E_k(\lambda_{Q_t}) ~dx
			\nonumber
			\\
			\geq&\frac{1}{k+1}V_k(Q_t)\log\frac{V_k(L)}{V_k(Q_t)}
			+ N_2^{k-1}\left(\frac{1}{2}(k-1+\frac{2k}{n})-\eta\right)\delta_2^2(N_2 L,Q_t).
		\end{align}
		Here $N_1=\left(\frac{V_k(Q_t)}{V_k(K)}\right)^{\frac{1}{k+1}}=\left(\frac{V_k(Q_t)}{V_k(L)}\right)^{\frac{1}{k+1}}=N_2$, since $V_k(K)=V_k(L)$. 	
		Multiplying equation \eqref{equ:pf-thm:p-BM-L2-Vk=Vk-4} by $(1-t)$
		and adding it to equation \eqref{equ:pf-thm:p-BM-L2-Vk=Vk-5} multiplied by $t$ yields,
		\begin{align}
			\label{equ:pf-thm:p-BM-L2-Vk=Vk-3}
			&- \frac{1}{k+1}V_k(Q_t)\log\frac{V_k(K)}{V_k(Q_t)} \nonumber
			\\
			\geq & N_1^{k-1}\left(\frac{1}{2}(k-1+\frac{2k}{n})-\eta\right)\left((1-t)\delta_2^2(N_1 K,Q_t)+t\delta_2^2(N_1 L,Q_t)\right).
		\end{align}

		On the other hands, applying the triangle inequality and H\"older's inequality, we obtain
		\begin{align*}
			\delta_2^2(K,L)\leq& \left(\delta_2(K,Q_t)+\delta_2(Q_t,L)\right)^2
			\nonumber
			\\
			=&\delta_2^2(K,Q_t)+\delta_2^2(Q_t,L)+2\delta_2(K,Q_t)\cdot \delta_2(Q_t,L)
			\nonumber
			\\
			\leq &(1+\frac{1-t}{t})\delta_2^2(K,Q_t)+(1+\frac{t}{1-t})\delta_2^2(Q_t,L)
			\nonumber
			\\
			=&
			\frac{1}{t(1-t)}\left((1-t)\delta_2^2(K,Q_t)+t\delta_2^2(Q_t,L)\right),
		\end{align*}
		where $t\in(0,1)$.  This means for $t\in[0,1]$, we have
		\begin{align}
			&\label{equ:pf-thm:p-BM-L2-Vk=Vk-2}
			(1-t)\delta_2^2(N_1K,Q_t)+t\delta_2^2(Q_t,N_1L)\nonumber
			\\
			\geq& t(1-t)\delta_2^2(N_1K,N_1L)\nonumber
			\\
			=&N_1^2t(1-t)\delta_2^2(K,L).
		\end{align}
		Putting \eqref{equ:pf-thm:p-BM-L2-Vk=Vk-2} into \eqref{equ:pf-thm:p-BM-L2-Vk=Vk-3} yields \eqref{equ:p-BM-Vk=p=0}.
		
		Next, we consider the case $p>-k+1-\frac{2k}{n}$ and $p\neq 0$. From Lemma \ref{lemma:C2-1-p-sum-C2-1},  for  sufficiently small $\varepsilon$,  $K,L\in\mathcal{N}_{\varepsilon}$, and $t\in[0,1]$, we have $Q_t:=(1-t)\cdot K+_pt\cdot L\in\mathcal{N}_{O(\varepsilon)}$ are convex bodies with  support functions $u_{Q_t}=\left((1-t)u_K^p+tu_L^p\right)^{\frac{1}{p}}$. Since $\mathcal{S}_p(K)=\mathcal{S}_p(L)=o$, we have  $\mathcal{S}_p(Q_t)=o$.
		Then
		\begin{align}\label{equ:pf-thm:p-BM-L2-Vk=Vk-1}
			&\frac{1}{p}V_k(Q_t)=\frac{1}{p}\int_{\mathbb{S}^n}u_{Q_t}E_k(\lambda_{Q_t})~dx\nonumber
			\\
			=&\frac{1-t}{p}\int_{\mathbb{S}^n}u_{K}^pu_{Q_t}^{1-p}E_k(\lambda_{Q_t})~dx
			+\frac{t}{p}\int_{\mathbb{S}^n}u_{L}^pu_{Q_t}^{1-p}E_k(\lambda_{Q_t})~dx\nonumber
			\\
			\overset{\eqref{equ:p-M-L2-p=!1}}{\geq}
			&\frac{1-t}{p}V_k(Q_t)^{\frac{k+1-p}{k+1}}V_k(K)^{\frac{p}{k+1}}+\frac{t}{p}V_k(Q_t)^{\frac{k+1-p}{k+1}}V_k(L)^{\frac{p}{k+1}}
			\nonumber
			\\
			&+{(1-t)N_1^{k-1-p}}\left(\frac{1}{2}(p+k-1+\frac{2k}{n})+O(\varepsilon)\right)\delta_2^2(N_1K,Q_t)\nonumber
			\\
			&+{tN_2^{k-1-p}}\left(\frac{1}{2}(p+k-1+\frac{2k}{n})+O(\varepsilon)\right)\delta_2^2(Q_t,N_2L)\nonumber
			\\
			=&\frac{1}{p}V_k(Q_t)^{\frac{k+1-p}{k+1}}\left((1-t)V_k(K)^{\frac{p}{k+1}}+tV_k(L)^{\frac{p}{k+1}}\right)
			\nonumber
			\\
			&+{N_1^{k-1-p}}\left(\frac{1}{2}(p+k-1+\frac{2k}{n})+O(\varepsilon)\right)\left((1-t)\delta_2^2(N_1K,Q_t)+t\delta_2^2(Q_t,N_1L)\right).
		\end{align}
		Here $N_1=\left(\frac{V_k(Q_t)}{V_k(K)}\right)^{\frac{1}{k+1}}=\left(\frac{V_k(Q_t)}{V_k(L)}\right)^{\frac{1}{k+1}}=N_2$, since $V_k(K)=V_k(L)$.
		
		Putting \eqref{equ:pf-thm:p-BM-L2-Vk=Vk-2} into \eqref{equ:pf-thm:p-BM-L2-Vk=Vk-1} yields \eqref{equ:p-BM-Vk=p!=0}.
	\end{proof}
	
	The unnormalized version of Theorem \ref{thm:p-BM-L2-Vk=Vk} is the following Theorem \ref{thm:p-BM-L2-Vk!=Vk}.
	\begin{theorem}
		\label{thm:p-BM-L2-Vk!=Vk}
		Let $K, L \in \mathcal{K}_o^{n+1}$ be  $C^3$  convex bodies, $1\leq k\leq n$, and $t\in[0,1]$. 
		For all $\eta>0$, there exists $\varepsilon>0$  such that if 
		$K,L\in\mathcal{N}_{\varepsilon}$ and $K,L$ satisfy  $\mathcal{S}_p(K)=\mathcal{S}_p(L)=o$, then we have for $p=0$,
		\begin{align}\label{equ:p-BM-Vk!=p=0}
			&\log	\left(V_k\((1-t)\cdot K+_0t\cdot L\)\right)
			-\log\( V_k(K)^{1-t} V_k(L)^{t}\)\nonumber
			\\
			\geq&
			(k+1)t(1-t)|\mathbb{S}^{n}|^{-1}\left(\frac{1}{2}(k-1+\frac{2k}{n})-\eta\right)\delta_2^2(N_K K,N_L L);
		\end{align}
		and for  $p
		\in(-k+1-\frac{2k}{n},0)\cup(0,+\infty)$, 
		\begin{align}\label{equ:p-BM-Vk!=p!=0}
			&\frac{1}{p}\left(V_k\((1-t)\cdot K+_pt\cdot L\)^{\frac{p}{k+1}}-\left((1-t)V_k(K)^{\frac{p}{k+1}}+tV_k(L)^{\frac{p}{k+1}}\right)\right)\nonumber
			\\
			\geq&t(1-t)c_0\cdot \left(\frac{1}{2}(p+k-1+\frac{2k}{n})-\eta\right)\delta_2^2(N_KK,N_LL).
		\end{align}	
		Here positive constants $N_L=\left(\frac{|\mathbb{S}^{n}|}{V_k(L)}\right)^{\frac{1}{k+1}}, $ $N_K=\left(\frac{|\mathbb{S}^{n}|}{V_k(K)}\right)^{\frac{1}{k+1}}$, and $$c_0=V_k(K)^{\frac{p}{k+1}}V_k(L)^{\frac{p}{k+1}} \left((1-t)V_k(K)^{\frac{p}{k+1}}+tV_k(L)^{\frac{p}{k+1}}\right)^{-1}|\mathbb{S}^{n}|^{-1}.$$
	\end{theorem}
	
	\begin{proof}
		We first consider the case $p=0$. If we take $\tilde{K}=N_K  K,\tilde{L}=N_L  L$, where $N_K=\left(\frac{V_k(\mathbb{B}^{n+1})}{V_k(K)}\right)^{\frac{1}{k+1}}$ and $N_L=\left(\frac{V_k(\mathbb{B}^{n+1})}{V_k(L)}\right)^{\frac{1}{k+1}}$, then $V_k(\tilde{K})=V_k(\tilde{L})=V_k(\mathbb{B}^{n+1})=|\mathbb{S}^n|$ and $\tilde{K},\tilde{L}$ are also  $C^3$  convex bodies in $\mathcal{N}_{O(\varepsilon)}$ by Lemma \ref{lemma:C2-1-p-sum-C2-1}, and one can check that $\tilde{K}$ and $\tilde{L}$ satisfy $\mathcal{S}_p(\tilde{K})=\mathcal{S}_p(\tilde{L})=o$. 
		Take
		\begin{align*}
			\beta=N_K^{1-t}N_L^t,
		\end{align*}
		then $(1-t)\cdot \tilde{K}+_0t\cdot \tilde{L}=\beta \left((1-t)\cdot K+_0t\cdot L\right)$  by Lemma \ref{lemma:C2-1-p-sum-C2-1} (iii).
		Form \eqref{equ:p-BM-Vk=p=0}, we have
		\begin{align*}
			&\log	\left(V_k\((1-t)\cdot \tilde{K}+_0t\cdot \tilde{L}\)\right)-\log\left(V_k(\tilde{K})\right)\nonumber
			\\
			\geq&
			(k+1)t(1-t)V_k(\tilde{K})^{-1}\left(\frac{1}{2}(k-1+\frac{2k}{n})-\eta\right)\delta_2^2(\tilde{K},\tilde{L}).
		\end{align*}
		This means
		\begin{align*}
			&\log	\left(V_k\((1-t)\cdot K+_0t\cdot L\)\right)
			+(1-t)\left(\log\(V_k(\mathbb{B}^{n+1})\)-\log\( V_k(K)\)\right)
			\\
			&+t\left(\log\(V_k(\mathbb{B}^{n+1})\)-\log\( V_k(L)\)\right)
			-\log\(V_k(\mathbb{B}^{n+1})\)\nonumber
			\\
			\geq&			(k+1)t(1-t)V_k(\mathbb{B}^{n+1})^{-1}\left(\frac{1}{2}(k-1+\frac{2k}{n})-\eta\right)\delta_2^2(\tilde{K},\tilde{L}),
		\end{align*}
		which implies \eqref{equ:p-BM-Vk!=p=0}.

		When $p\neq 0$, we take
		\begin{align*}
			\tilde{t}=\frac{tV_k(L)^{\frac{p}{k+1}}}{(1-t)V_k(K)^{\frac{p}{k+1}}+tV_k(L)^{\frac{p}{k+1}}}\in[0,1].
		\end{align*}
		Form \eqref{equ:p-BM-Vk=p!=0}, 
		 we have
		\begin{align*}
			&\frac{1}{p}\left(V_k\((1-\tilde{t})\cdot \tilde{K}+_p\tilde{t}\cdot \tilde{L}\)^{\frac{p}{k+1}}-V_k(\tilde{L})^{\frac{p}{k+1}}\right)\nonumber
			\\
			\geq&\tilde{t}(1-\tilde{t})V_k(\tilde{K})^{-\frac{k+1-p}{k+1}}\left(\frac{1}{2}(p+k-1+\frac{2k}{n})-\eta\right)\delta_2^2(\tilde{K},\tilde{L}).
		\end{align*}
		This means
		\begin{align*}
			&\frac{1}{p}\left(\frac{V_k\((1-t)\cdot K+_pt\cdot L\)^{\frac{p}{k+1}}V_k(\mathbb{B}^{n+1})^{\frac{p}{k+1}}}{(1-t)V_k(K)^{\frac{p}{k+1}}+tV_k(L)^{\frac{p}{k+1}}}-V_k(\mathbb{B}^{n+1})^{\frac{p}{k+1}}\right)\nonumber
			\\
			\geq& \frac{t(1-t)V_k(K)^{\frac{p}{k+1}}V_k(L)^{\frac{p}{k+1}}}{\left((1-t)V_k(K)^{\frac{p}{k+1}}+tV_k(L)^{\frac{p}{k+1}}\right)^2}V_k(\mathbb{B}^{n+1})^{\frac{p}{k+1}-1}\left(\frac{1}{2}(p+k-1+\frac{2k}{n})-\eta\right)\delta_2^2(\tilde{K},\tilde{L}),
		\end{align*}
		which proves \eqref{equ:p-BM-Vk!=p!=0}. 	
	\end{proof}

	From Corollary \ref{thm:p-M-index-k!=n=!Vk}, we can derive Theorem \ref{thm:p-BM-L2-Vk=Vk-index}.	
	\begin{proof}[\textbf{Proof of Theorem \ref{thm:p-BM-L2-Vk=Vk-index}}]
	  We first consider the case $p=0$. From Lemma \ref{lemma:C2-1-p-sum-C2-1},  for  sufficiently small $\varepsilon$,  $K,L\in\mathcal{N}_{\varepsilon}$, and $t\in[0,1]$, we have $Q_t:=(1-t)\cdot K+_0t\cdot L\in\mathcal{N}_{O(\varepsilon)}$ is a convex body with  support function $u_{Q_t}=u_K^{1-t} u_L^{t}$ and $\mathcal{S}_p(Q_t)=o$.   By \eqref{equ:p-M-index-k!=n-p=0!=Vk}, we have
		\begin{align}
			\label{equ:pf-thm:p-BM-L2-Vk=Vk-4-index}
			&-t\int_{\mathbb{S}^n}\left(\log(\frac{u_L}{u_{K}})\right)u_{Q_t}E_k(\lambda_{Q_t}) ~dx=\int_{\mathbb{S}^n}\left(\log(\frac{u_K}{u_{Q_t}})\right)u_{{Q_t}}E_k(\lambda_{Q_t}) ~dx
			\nonumber
			\\
			\geq&\frac{1}{k+1}V_k(Q_t)\log\frac{V_k(K)}{V_k(Q_t)}
			+  V_k(Q_t)\left(\frac{(k-1)n+2k}{2n(n+1)^2}-\eta\right)\alpha(K,Q_t)^2,
		\end{align} 	
		and
		\begin{align}
			\label{equ:pf-thm:p-BM-L2-Vk=Vk-5-index}
			&(1-t)\int_{\mathbb{S}^n}\left(\log(\frac{u_L}{u_K})\right)u_{Q_t}E_k(\lambda_{Q_t}) ~dx=\int_{\mathbb{S}^n}\left(\log(\frac{u_L}{u_{Q_t}})\right)u_{{Q_t}}E_k(\lambda_{Q_t}) ~dx
			\nonumber
			\\
			\geq&\frac{1}{k+1}V_k(Q_t)\log\frac{V_k(L)}{V_k(Q_t)}
			+  V_k(Q_t)\left(\frac{(k-1)n+2k}{2n(n+1)^2}-\eta\right)\alpha(L,Q_t)^2.
		\end{align}
		Multiplying equation \eqref{equ:pf-thm:p-BM-L2-Vk=Vk-4-index} by $(1-t)$
		and adding it to equation \eqref{equ:pf-thm:p-BM-L2-Vk=Vk-5-index} multiplied by $t$ yields,
		\begin{align}
			\label{equ:pf-thm:p-BM-L2-Vk=Vk-3-index}
			&- \frac{1}{k+1}V_k(Q_t)\log\frac{V_k(K)}{V_k(Q_t)} \nonumber
			\\
			\geq &  V_k(Q_t)\left(\frac{(k-1)n+2k}{2n(n+1)^2}-\eta\right)\left((1-t)\alpha(K,Q_t)^2+t\alpha(L,Q_t)^2\right).
		\end{align}	
		
		On the other hands, applying the triangle inequality of $\alpha(K,L)$ \cite{Figalli-Maggi-Pratelli-2010}*{eq. (4.2)} and H\"older's inequality, we obtain
		\begin{align*}
			\alpha(K,L)^2\leq& \left(\alpha(K,Q_t)+\alpha(Q_t,L)\right)^2
			\nonumber
			\\
			=&\alpha(K,Q_t)^2+\alpha(Q_t,L)^2+2\alpha(K,Q_t)\cdot \alpha(Q_t,L)
			\nonumber
			\\
			\leq &
			\frac{1}{t(1-t)}\left((1-t)\alpha(K,Q_t)^2+t\alpha(Q_t,L)^2\right),
		\end{align*}
		where $t\in(0,1)$.  This means that for $t\in[0,1]$, we have
		\begin{align}
			\label{equ:pf-thm:p-BM-L2-Vk=Vk-2-index}
			(1-t)\alpha(K,Q_t)^2+t\alpha(L,Q_t)^2
			\geq t(1-t)\alpha(K,L)^2.
		\end{align}
		Putting \eqref{equ:pf-thm:p-BM-L2-Vk=Vk-2-index} into \eqref{equ:pf-thm:p-BM-L2-Vk=Vk-3-index} yields \eqref{equ:p-BM-Vk=p=0-index}.
		
		Next, we consider the case $p>-k+1-\frac{2k}{n}$ and $p\neq 0$. From Lemma \ref{lemma:C2-1-p-sum-C2-1},  for  sufficiently small $\varepsilon$,  $K,L\in\mathcal{N}_{\varepsilon}$, and $t\in[0,1]$, we have $Q_t:=(1-t)\cdot K+_pt\cdot L\in\mathcal{N}_{O(\varepsilon)}$ is a convex body with  support function $u_{Q_t}=\left((1-t)u_K^p+tu_L^p\right)^{\frac{1}{p}}$, and
		$\mathcal{S}_p(Q_t)=o$.
		Then
		\begin{align}\label{equ:pf-thm:p-BM-L2-Vk=Vk-1-index}
			&\frac{1}{p}V_k(Q_t)=\frac{1}{p}\int_{\mathbb{S}^n}u_{Q_t}E_k(\lambda_{Q_t})~dx\nonumber
			\\
			=&\frac{1-t}{p}\int_{\mathbb{S}^n}u_{K}^pu_{Q_t}^{1-p}E_k(\lambda_{Q_t})~dx
			+\frac{t}{p}\int_{\mathbb{S}^n}u_{L}^pu_{Q_t}^{1-p}E_k(\lambda_{Q_t})~dx\nonumber
			\\
			\overset{\eqref{equ:p-M-index-k!=n-p!=0!=Vk}}{\geq}
			&\frac{1-t}{p}V_k(Q_t)^{\frac{k+1-p}{k+1}}V_k(K)^{\frac{p}{k+1}}+\frac{t}{p}V_k(Q_t)^{\frac{k+1-p}{k+1}}V_k(L)^{\frac{p}{k+1}}
			\nonumber
			\\
			&+(1-t)N_1^{-p}V_k(Q_t)\left(\frac{n(k-1+p)+2k}{2n(n+1)^2}-{\eta}\right)\alpha(K,Q_t)^2\nonumber
			\\
			&tN_2^{-p}V_k(Q_t)\left(\frac{n(k-1+p)+2k}{2n(n+1)^2}-{\eta}\right)\alpha(L,Q_t)^2\nonumber
			\\
			=&\frac{1}{p}V_k(Q_t)^{\frac{k+1-p}{k+1}}\left((1-t)V_k(K)^{\frac{p}{k+1}}+tV_k(L)^{\frac{p}{k+1}}\right)
			\nonumber
			\\
			&+N_1^{-p}V_k(Q_t)\left(\frac{n(k-1+p)+2k}{2n(n+1)^2}-{\eta}\right)\left((1-t)\alpha(K,Q_t)^2+t\alpha(K,Q_t)^2\right),
		\end{align}
		where  $N_1=\left(\frac{V_k(Q_t)}{V_k(K)}\right)^{\frac{1}{k+1}}=\left(\frac{V_k(Q_t)}{V_k(L)}\right)^{\frac{1}{k+1}}=N_2$, since $V_k(K)=V_k(L)$. Then $N_1^{-p}V_k(Q_t)=V_k(Q_t)^{1-\frac{p}{k+1}}V_{k}(K)^{\frac{p}{k+1}}$.	Putting \eqref{equ:pf-thm:p-BM-L2-Vk=Vk-2-index} into \eqref{equ:pf-thm:p-BM-L2-Vk=Vk-1-index} yields \eqref{equ:p-BM-Vk=p!=0-index}.

	\end{proof}
	
	Similar to Theorem \ref{thm:p-BM-L2-Vk!=Vk},  Theorem \ref{thm:p-BM-L2-Vk=Vk-index}   implies the following theorems.
	\begin{theorem}
		\label{thm:p-BM-L2-Vk!=Vk-index}
		Let $K, L \in \mathcal{K}_o^{n+1}$ be  $C^3$  convex bodies, $1\leq k\leq n$, and $t\in[0,1]$.  For all $\eta>0$, there exists $\varepsilon>0$  such that if $K,L\in\mathcal{N}_{\varepsilon}$ and $K,L$ satisfy  $\mathcal{S}_p(K)=\mathcal{S}_p(L)=o$, then we have for $p=0$,
		\begin{align}\label{equ:p-BM-Vk!=p=0-index}
			&\log	\left(V_k\((1-t)\cdot K+_0t\cdot L\)\right)
			-\log\( V_k(K)^{1-t} V_k(L)^{t}\)\nonumber
			\\
			\geq&
			(k+1)t(1-t)\left(\frac{(k-1)n+2k}{2n(n+1)^2}-\eta\right)\alpha(K,L)^2;
		\end{align}
	and for $p\in(1-k-\frac{2k}{n},0)\cup(0, +\infty)$, 
		\begin{align}\label{equ:p-BM-Vk!=p!=0-index}
			&\frac{1}{p}\left(V_k\((1-t)\cdot K+_pt\cdot L\)^{\frac{p}{k+1}}-(1-t)V_k(K)^{\frac{p}{k+1}}-tV_k(L)^{\frac{p}{k+1}}\right)\nonumber
			\\
			\geq&t(1-t)c_1\left(\frac{n(k-1+p)+2k}{2n(n+1)^2}-{\eta}\right)\alpha(K,L)^2.
		\end{align}
		Here
		\begin{align*}
			c_1=&
			V_k(K)^{\frac{p}{k+1}}
			V_k(L)^{\frac{p}{k+1}}\left((1-t)V_k(K)^{\frac{p}{k+1}}+tV_k(L)^{\frac{p}{k+1}}\right)^{-1}.
		\end{align*}
	\end{theorem}
	
	\section{A uniqueness result for the $L^p$-Christoffel-Minkowski problem}\label{sec 6}
	The stability results for the inequalities imply the corresponding inequalities themselves. Using the $L^p$-Minkowski (or $L^p$-Brunn--Minkowski) inequalities for intrinsic volumes, we then derive a uniqueness result for the $L^p$-Christoffel-Minkowski problem.

	\begin{proof}[\textbf{Proof of Corollary \ref{thm:Christoffel} (c)}]
		 When $p=0$,  from Corollary \ref{thm:Christoffel} (a), 
		we obtain
		\begin{align}
			\int_{\mathbb{S}^n} \Bigl(\log \frac{u_L}{u_K}\Bigr) u_K E_k(\lambda_K) ~dx
			- \frac{1}{k+1}V_k(K) \log\frac{V_k(L)}{V_k(K)} \geq 0. \label{equ:pf-Chri-2}
		\end{align}
		Interchanging $K$ and $L$ in Corollary \ref{thm:Christoffel} (a) gives
		\begin{align*}
			\int_{\mathbb{S}^n} \Bigl(\log \frac{u_K}{u_L}\Bigr) u_L E_k(\lambda_L) ~ dx
			- \frac{1}{k+1}V_k(L) \log\frac{V_k(K)}{V_k(L)} \geq  0. 
		\end{align*}
		Using the condition \eqref{equ:Christolell-equ} for $p=0$, we derive that
		\begin{align}
			-\int_{\mathbb{S}^n} \Bigl(\log \frac{u_L}{u_K}\Bigr) u_K E_k(\lambda_K) ~dx
			- \frac{1}{k+1}V_k(L) \log\frac{V_k(K)}{V_k(L)} \geq  0, \label{equ:pf-Chri-4}
		\end{align}
		and
		\begin{align}
			V_k(K)=\int_{\mathbb{S}^n}u_KE_k(\lambda_K)~dx\overset{\eqref{equ:Christolell-equ}}{=}\int_{\mathbb{S}^n}u_LE_k(\lambda_L)~dx= V_k(L). \label{equ:pf-Chri-7}
		\end{align}
		Substituting \eqref{equ:pf-Chri-7} into \eqref{equ:pf-Chri-2} gives
		\begin{align}
			\int_{\mathbb{S}^n} \Bigl(\log \frac{u_L}{u_K}\Bigr) u_K E_k(\lambda_K) ~dx \geq  0, \label{equ:pf-Chri-8}
		\end{align}
		while substituting it into \eqref{equ:pf-Chri-4} gives
		\begin{align*}
			-\int_{\mathbb{S}^n} \Bigl(\log \frac{u_L}{u_K}\Bigr) u_K E_k(\lambda_K) ~dx \geq  0. 
		\end{align*}
		Hence the integral must be zero, and equality holds in \eqref{equ:pf-Chri-2}.
		By Corollary \ref{thm:Christoffel} (a), equality in \eqref{equ:pf-Chri-2} occurs if and only if $K = cL$. Putting it into \eqref{equ:Christolell-equ} we have $c=1$.
		\\
		
		 When $p\neq 0$, by \eqref{equ:Christolell-equ}, we can derive
		\begin{align}\label{equ:pf-christ-0}
			\int_{\mathbb{S}^n}u_L^pu_K^{1-p}E_k(\lambda_K)~dx=\int_{\mathbb{S}^n}u_L^pu_L^{1-p}E_k(\lambda_L)~dx=V_k(L).
		\end{align}
		Substituting \eqref{equ:pf-christ-0} into \eqref{equ:p-M-inequ-p!=0}, we obtain equality when $p= k+1$, and the equality case implies $K=cL$.
		
		If $p\neq k+1$, putting \eqref{equ:pf-christ-0} into \eqref{equ:p-M-inequ-p!=0}, we have
		\begin{align}\label{equ:pf-christ-1}
			\frac{1}{p}V_k(L)\left(1-\left(\frac{V_k(K)}{V_k(L)}\right)^{\frac{k+1-p}{k+1}}\right)\geq 0.
		\end{align}
		The same conclusion holds upon interchanging $K$ and $L$  in Corollary \ref{thm:Christoffel} (a). Analogously, we have
		\begin{align}
			\label{equ:pf-christ-2}		\frac{1}{p}V_k(K)\left(1-\left(\frac{V_k(L)}{V_k(K)}\right)^{\frac{k+1-p}{k+1}}\right)\geq 0.
		\end{align}
		Combining \eqref{equ:pf-christ-1} and \eqref{equ:pf-christ-2}, we obtain
		\begin{align*}
			V_k(K)=V_k(L).
		\end{align*}
		Then
		\begin{align}\label{equ:pf-christ-3}
			V_k(K)^{\frac{k+1-p}{k+1}}V_k(L)^{\frac{p}{k+1}}	=V_k(L).
		\end{align}
		Observing \eqref{equ:pf-christ-0} and \eqref{equ:pf-christ-3} shows that equality holds in \eqref{equ:p-M-inequ-p!=0}, and hence by Corollary \ref{thm:Christoffel} (a), $K=cL$. Substituting this into \eqref{equ:Christolell-equ} yields  $c=1$ for $1-p+k\neq 0$.
	\end{proof}

	\vskip 2mm  {\bf Acknowledgment.} 
	The research of  JG was supported by NSF of Zhejiang Province No. LQN26A010008.
	The research of  GL was supported by NSFC No. 12261105.
	
	\vskip 2mm
	\textbf{Data Availability.} Data sharing is not applicable to this article as no datasets were generated or analysed during the current study.
	
	\vskip 2mm
	\textbf{Conflict of interest.} On behalf of all authors, the corresponding author states that there is no conflict of interest.


	\section*{Reference}
	

\bigbreak

\noindent Shanwei Ding,  {\tt dingsw@whu.edu.cn} \\
School of Mathematics and Statistics, Huazhong University of Science and
Technology, Wuhan, Hubei 430074, China\\

\noindent  Jinyu Gao, {\tt gaojy@zstu.edu.cn} \\
Department of Mathematics,
Zhejiang Sci-Tech University,
Hangzhou, Zhejiang 310018, China\\

\noindent  Guanghan Li, {\tt ghli@whu.edu.cn} \\
School of Mathematics and Statistics, Wuhan University,
Wuhan, Hubei 430072, China


\begin{biblist}
		
		
		
		\bib{Andrews-97}{article}{
			author={Andrews, B.},
			title={Monotone quantities and unique limits for evolving convex
				hypersurfaces},
			journal={Internat. Math. Res. Notices},
			date={1997},
			number={20},
			pages={1001--1031},
			issn={1073-7928},
			review={\MR{1486693}},
			doi={10.1155/S1073792897000640},
		}
		\bib{Andrews-Guan-Ni-2016}{article}{
			author={Andrews, B.},
			author={Guan, P.},
			author={Ni, L.},
			title={Flow by powers of the Gauss curvature},
			journal={Adv. Math.},
			volume={299},
			date={2016},
			pages={174--201},
			issn={0001-8708},
			review={\MR{3519467}},
			doi={10.1016/j.aim.2016.05.008},
		}
		
		\bib{Bianchini-Colesanti-Pagnini-Roncoroni-Math-ann-2023}{article}{
			author={Bianchini, C.},
			author={Colesanti, A.},
			author={Pagnini, D.},
			author={Roncoroni, A.},
			title={On $p$-Brunn--Minkowski inequalities for intrinsic volumes, with
				$0\leq p < 1$},
			journal={Math. Ann.},
			volume={387},
			date={2023},
			number={1-2},
			pages={321--352},
			issn={0025-5831},
			review={\MR{4631049}},
			doi={10.1007/s00208-022-02454-0},
		}
		
		
		\bib{Borell-1997}{article}{
			author={Borell, C.},
			title={The Brunn--Minkowski inequality in Gauss space},
			journal={Invent. Math.},
			volume={30},
			date={1975},
			number={2},
			pages={207--216},
			issn={0020-9910},
			review={\MR{399402}},
			doi={10.1007/BF01425510},
		}
		
		
		\bib{Boroczky-Chen-Liu-Saroglou-arxiv}{article}{
			title={Uniqueness in the near isotropic $L_p$ dual Minkowski problem},
			author={B\"{o}r\"{o}czky, K. J.},
			author={S. Chen},
			author={W. Liu},
			author={C. Saroglou},
			year={2025},
			eprint={2505.01066},
			archivePrefix={arXiv},
			primaryClass={math.AP},
			url={https://arxiv.org/abs/2505.01066},
		}
		\bib{Boroczky-De-2024}{article}{
			author={B\"{o}r\"{o}czky, K. J.},
			author={De, A.},
			title={Stability of the logarithmic Brunn--Minkowski inequality in the
				case of many hyperplane symmetries},
			journal={Studia Sci. Math. Hungar.},
			volume={61},
			date={2024},
			number={4},
			pages={341--376},
			issn={0081-6906},
			review={\MR{4847516}},
		}
		\bib{Boroczky-Figalli-Ramos-2025-book}{book}{%
			author={B\"{o}r\"{o}czky, K. J.},
			author={Figalli, A.},
			author={Ramos, J. P. G.},
			title={Isoperimetric inequalities, Brunn-Minkowski theory and
				Minkowski-type Monge-Amp\`ere equations on the sphere},
			series={Zurich Lectures in Advanced Mathematics},
			publisher={EMS Press, Berlin},
			date={[2026] \copyright 2026},
			pages={xv+523},
			isbn={978-3-98547-104-1},
			isbn={978-3-98547-604-6},
			review={\MR{5095086}},
		}
		\bib{Boroczky-Kalantzopoulos-2022}{article}{
			author={B\"{o}r\"{o}czky, K. J.},
			author={Kalantzopoulos, P.},
			title={Log-Brunn--Minkowski inequality under symmetry},
			journal={Trans. Amer. Math. Soc.},
			volume={375},
			date={2022},
			number={8},
			pages={5987--6013},
			issn={0002-9947},
			review={\MR{4469244}},
			doi={10.1090/tran/8691},
		}
		
		\bib{Boroczky-Lutwak-Yang-Zhang-2012}{article}{
			author={B\"{o}r\"{o}czky, K. J.},
			author={Lutwak, E.},
			author={Yang, D.},
			author={Zhang, G.},
			title={The log-Brunn--Minkowski inequality},
			journal={Adv. Math.},
			volume={231},
			date={2012},
			number={3-4},
			pages={1974--1997},
			issn={0001-8708},
			review={\MR{2964630}},
			doi={10.1016/j.aim.2012.07.015},
		}
		
		\bib{Boroczky-Lutwak-Yang-Zhang-2013}{article}{
			author={B\"{o}r\"{o}czky, K. J.},
			author={Lutwak, E.},
			author={Yang, D.},
			author={Zhang, G.},
			title={The logarithmic Minkowski problem},
			journal={J. Amer. Math. Soc.},
			volume={26},
			date={2013},
			number={3},
			pages={831--852},
			issn={0894-0347},
			review={\MR{3037788}},
			doi={10.1090/S0894-0347-2012-00741-3},
		}
		\bib{Boroczky-Machado-Ramos-2026-arxiv}{article}{%
			title={Quantitative Stability for Minkowski's problem},
			author={K. J. B{\"o}r{\"o}czky},
			author={J. Machado},
			author={J. P. G. Ramos},
			year={2026},
			eprint={2603.17726},
			archivePrefix={arXiv},
			primaryClass={math.AP},
			url={https://arxiv.org/abs/2603.17726},
		}
		
		\bib{Boroczky-Saroglou-Cvpde-2024}{article}{
			author={B\"{o}r\"{o}czky, K. J.},
			author={Saroglou, C.},
			title={Uniqueness when the $L_p$ curvature is close to be a constant for
				$p\in [0,1)$},
			journal={Calc. Var. Partial Differential Equations},
			volume={63},
			date={2024},
			number={6},
			pages={Paper No. 154, 26},
			issn={0944-2669},
			review={\MR{4765819}},
			doi={10.1007/s00526-024-02763-z},
		}
		\bib{Brendle-Choi-Daskalopoulos}{article}{
			author={Brendle, S.},
			author={Choi, K.},
			author={Daskalopoulos, P.},
			title={Asymptotic behavior of flows by powers of the Gaussian curvature},
			journal={Acta Math.},
			volume={219},
			date={2017},
			number={1},
			pages={1--16},
			issn={0001-5962},
			review={\MR{3765656}},
			doi={10.4310/ACTA.2017.v219.n1.a1},
		}
\bib{Chen-Li-2021}{article}{
			author={Chen, H.},
			author={Li, Q.},
			title={The $L_ p$ dual Minkowski problem and related parabolic flows},
			journal={J. Funct. Anal.},
			volume={281},
			date={2021},
			number={8},
			pages={Paper No. 109139, 65},
			issn={0022-1236},
			review={\MR{4271790}},
			doi={10.1016/j.jfa.2021.109139},
		}
		\bib{Chen-JFA-2020}{article}{
			author={Chen, L.},
			title={Uniqueness of solutions to $L_p$-Christoffel-Minkowski problem for
				$p<1$},
			journal={J. Funct. Anal.},
			volume={279},
			date={2020},
			number={8},
			pages={108692, 15},
			issn={0022-1236},
			review={\MR{4116152}},
			doi={10.1016/j.jfa.2020.108692},
		}
\bib{Chen-Huang-Li-Liu-2020}{article}{
			author={Chen, S.},
			author={Huang, Y.},
			author={Li, Q.-R.},
			author={Liu, J.},
			title={The $L_ p$-Brunn--Minkowski inequality for $p<1$},
			journal={Adv. Math.},
			volume={368},
			date={2020},
			pages={107166, 21},
			issn={0001-8708},
			review={\MR{4088419}},
			doi={10.1016/j.aim.2020.107166},
		}
		\bib{Cicalese-Leonardi-2012}{article}{
			author={Cicalese, M.},
			author={Leonardi, G.},
			title={A selection principle for the sharp quantitative isoperimetric
				inequality},
			journal={Arch. Ration. Mech. Anal.},
			volume={206},
			date={2012},
			number={2},
			pages={617--643},
			issn={0003-9527},
			review={\MR{2980529}},
			doi={10.1007/s00205-012-0544-1},
		}
			\bib{Colesanti-Livshyts-Livshyts-2017}{article}{
				author={Colesanti, A.},
				author={Livshyts, G. V.},
				author={Marsiglietti, A.},
				title={On the stability of Brunn--Minkowski type inequalities},
				journal={J. Funct. Anal.},
				volume={273},
				date={2017},
				number={3},
				pages={1120--1139},
				issn={0022-1236},
				review={\MR{3653949}},
				doi={10.1016/j.jfa.2017.04.008},
			}
			\bib{Ding-Li-Trans}{article}{
				author={Ding, S.},
				author={Li, G.},
				title={A class of inverse curvature flows and $L^p$ dual
					Christoffel-Minkowski problem},
				journal={Trans. Amer. Math. Soc.},
				volume={376},
				date={2023},
				number={1},
				pages={697--752},
				issn={0002-9947},
				review={\MR{4510121}},
				doi={10.1090/tran/8793},
			}
		\bib{Ding-Li-2025-JDE}{article}{
			author={Ding, S.},
			author={Li, G.},
			title={A type of anisotropic flows and dual Orlicz Christoffel-Minkowski
				type equations},
			journal={J. Differential Equations},
			volume={445},
			date={2025},
			pages={Paper No. 113586, 34},
			issn={0022-0396},
			review={\MR{4925888}},
			doi={10.1016/j.jde.2025.113586},
		}
		
		\bib{Diskant-1973}{article}{
			author={Diskant, V. I.},
			title={Stability of the solution of a Minkowski equation},
			language={Russian},
			journal={Siberian J. Math.},
			volume={14},
			date={1973},
			pages={669--673, 696},
			issn={0037-4474},
		}
		
		\bib{Figalli-Hintum-Tiba-arxiv-2023}{article}{
			title={Sharp quantitative stability of the Brunn--Minkowski inequality},
			author={A. Figalli},
			author={P. van  Hintum},
			author={M. Tiba},
			year={2023},
			eprint={2310.20643},
			archivePrefix={arXiv},
			primaryClass={math.AP},
			url={https://arxiv.org/abs/2310.20643},
		}
		
		\bib{Figalli-Hintum-Tiba-arxiv-2024}{article}{
			title={Sharp stability of the Brunn--Minkowski inequality via optimal mass transportation},
			author={A. Figalli},
			author={P. van  Hintum},
			author={M. Tiba},
			year={2024},
			eprint={2407.10932},
			archivePrefix={arXiv},
			primaryClass={math.AP},
			url={https://arxiv.org/abs/2407.10932},
		}
		
		\bib{Figalli-Maggi-Pratelli-2009}{article}{
			author={Figalli, A.},
			author={Maggi, F.},
			author={Pratelli, A.},
			title={A refined Brunn--Minkowski inequality for convex sets},
			journal={Ann. Inst. H. Poincar\'{e} C Anal. Non Lin\'{e}aire},
			volume={26},
			date={2009},
			number={6},
			pages={2511--2519},
			issn={0294-1449},
			review={\MR{2569906}},
			doi={10.1016/j.anihpc.2009.07.004},
		}
		
		
		\bib{Figalli-Maggi-Pratelli-2010}{article}{
			author={Figalli, A.},
			author={Maggi, F.},
			author={Pratelli, A.},
			title={A mass transportation approach to quantitative isoperimetric
				inequalities},
			journal={Invent. Math.},
			volume={182},
			date={2010},
			number={1},
			pages={167--211},
			issn={0020-9910},
			review={\MR{2672283}},
			doi={10.1007/s00222-010-0261-z},
		}
		\bib{Firey-book}{article}{
			author={Firey, W. J.},
			title={$p$-means of convex bodies},
			journal={Math. Scand.},
			volume={10},
			date={1962},
			pages={17--24},
			issn={0025-5521},
			review={\MR{141003}},
			doi={10.7146/math.scand.a-10510},
		}
		\bib{Fusco-Maggi-2008}{article}{
			author={Fusco, N.},
			author={Maggi, F.},
			author={Pratelli, A.},
			title={The sharp quantitative isoperimetric inequality},
			journal={Ann. of Math. (2)},
			volume={168},
			date={2008},
			number={3},
			pages={941--980},
			issn={0003-486X},
			review={\MR{2456887}},
			doi={10.4007/annals.2008.168.941},
		}
		\bib{Gao-Li-cvpde}{article}{
			author={Gao, J.},
			author={Li, G.},
			title={Alexandrov-Fenchel inequalities for convex anisotropic capillary
				hypersurfaces in the half-space},
			journal={Calc. Var. Partial Differential Equations},
			volume={65},
			date={2026},
			number={9},
			pages={Paper No. 254},
			issn={0944-2669},
			review={\MR{5122632}},
			doi={10.1007/s00526-026-03416-z},
		}
		\bib{Gardner-book-2006}{book}{
			author={Gardner, R. J.},
			title={Geometric tomography},
			series={Encyclopedia of Mathematics and its Applications},
			volume={58},
			edition={2},
			publisher={Cambridge University Press, New York},
			date={2006},
			pages={xxii+492},
			isbn={0-521},
			isbn={0-521-68493-5},
			review={\MR{2251886}},
			doi={10.1017/CBO9781107341029},
		}
		
		\bib{Gardner-2002}{article}{
			author={Gardner, R. J.},
			title={The Brunn--Minkowski inequality},
			journal={Bull. Amer. Math. Soc. (N.S.)},
			volume={39},
			date={2002},
			number={3},
			pages={355--405},
			issn={0273-0979},
			review={\MR{1898210}},
			doi={10.1090/S0273-0979-02-00941-2},
		}
		
		\bib{Groemer-1988}{article}{
			author={Groemer, H.},
			title={On the Brunn--Minkowski theorem},
			journal={Geom. Dedicata},
			volume={27},
			date={1988},
			number={3},
			pages={357--371},
			issn={0046-5755},
			review={\MR{960207}},
			doi={10.1007/BF00181500},
		}
		\bib{Groemer-1993}{article}{
			author={Groemer, H.},
			title={Stability of geometric inequalities},
			conference={
				title={Handbook of convex geometry, Vol. A, B},
			},
			book={
				publisher={North-Holland, Amsterdam},
			},
			isbn={0-444-89598-1},
			date={1993},
			pages={125--150},
		}
		\bib{Groemer-1991}{article}{
			author={Groemer, H.},
			author={Schneider, R.},
			title={Stability estimates for some geometric inequalities},
			journal={Bull. London Math. Soc.},
			volume={23},
			date={1991},
			number={1},
			pages={67--74},
			issn={0024-6093},
			review={\MR{1111537}},
			doi={10.1112/blms/23.1.67},
		}
		\bib{Guan-Ni-2017}{article}{
			author={Guan, P.},
			author={Ni, L.},
			title={Entropy and a convergence theorem for Gauss curvature flow in high
				dimension},
			journal={J. Eur. Math. Soc. (JEMS)},
			volume={19},
			date={2017},
			number={12},
			pages={3735--3761},
			issn={1435-9855},
			review={\MR{3730513}},
			doi={10.4171/JEMS/752},
		}
		 \bib{Xia-Guan-2018}{article}{
		 	author={Guan, P.},
		 	author={Xia, C.},
		 	title={$L^p$ Christoffel-Minkowski problem: the case $1<p<k+1$},
		 	journal={Calc. Var. Partial Differential Equations},
		 	volume={57},
		 	date={2018},
		 	number={2},
		 	pages={Paper No. 69, 23},
		 	issn={0944-2669},
		 	review={\MR{3776359}},
		 	doi={10.1007/s00526-018-1341-y},
		 }
		
		\bib{Ivaki-Milman-2023}{article}{
			author={M. N. Ivaki},
			author={E. Milman},
			title={ Uniqueness of solutions to a class of isotropic curvature problems},
			journal={Adv. Math.},
			volume={435},
			date={2023},
			pages={109350},
		}
		\bib{Kolesnikov-Milman-MAMS-2022}{article}{
			author={Kolesnikov, A. V.},
			author={Milman, E.},
			title={Local $L^p$-Brunn--Minkowski inequalities for $p<1$},
			journal={Mem. Amer. Math. Soc.},
			volume={277},
			date={2022},
			number={1360},
			pages={v+78},
			issn={0065-9266},
			isbn={978-1-4704-5160-8;978-1-4704-7092-0},
			review={\MR{4438690}},
			doi={10.1090/memo/1360},
		}
        \bib{Li-Wan-2024}{article}{
			author={Li, H.},
			author={Wan, Y.},
			title={Uniqueness of solutions to some classes of anisotropic and
				isotropic curvature problems},
			journal={J. Funct. Anal.},
			volume={287},
			date={2024},
			number={3},
			pages={Paper No. 110471, 30},
			issn={0022-1236},
			review={\MR{4742772}},
			doi={10.1016/j.jfa.2024.110471},
		}
		\bib{L93}{article}{
			author={Lutwak, E.},
			title={The Brunn--Minkowski-Firey theory. I. Mixed volumes and the
				Minkowski problem},
			journal={J. Differential Geom.},
			volume={38},
			date={1993},
			number={1},
			pages={131--150},
			issn={0022-040X},
			review={\MR{1231704}},
		}
		\bib{Lutwak-Yang-Zhang-2012}{article}{
			author={Lutwak, E.},
			author={Yang, D.},
			author={Zhang, G.},
			title={The Brunn--Minkowski-Firey inequality for nonconvex sets},
			journal={Adv. in Appl. Math.},
			volume={48},
			date={2012},
			number={2},
			pages={407--413},
			issn={0196-8858},
			review={\MR{2873885}},
			doi={10.1016/j.aam.2011.11.003},
		}
		
		\bib{McCoy-2011}{article}{
			author={McCoy, J. A.},
			title={Self-similar solutions of fully nonlinear curvature flows},
			journal={Ann. Sc. Norm. Super. Pisa Cl. Sci. (5)},
			volume={10},
			date={2011},
			number={2},
			pages={317--333},
			issn={0391-173X},
			review={\MR{2856150}},
		}
		\bib{Xia-arxiv}{article}{
			author={X. Mei},
			author={G. Wang},
			author={L. Weng},
			author={C. Xia},
			title={Alexandrov-Fenchel inequalities for convex hypersurfaces in the half-space with capillary boundary. II},
			journal={Math. Z.},
			volume={310},
			year={2025},
			number={4},
			pages={Paper No. 71.},
		}
		
		
		\bib{Milman-JDG-2024}{article}{
			author={Milman, E.},
			title={A sharp centro-affine isospectral inequality of
				Szeg\"{o}-Weinberger type and the $L^p$-Minkowski problem},
			journal={J. Differential Geom.},
			volume={127},
			date={2024},
			number={1},
			pages={373--408},
			issn={0022-040X},
			review={\MR{4753504}},
			doi={10.4310/jdg/1717356160},
		}
		
		\bib{Milman-JEMS-2025}{article}{
			author={Milman, E.},
			title={Centro-affine differential geometry and the log-Minkowski problem},
			journal={J. Eur. Math. Soc. (JEMS)},
			volume={27},
			date={2025},
			number={2},
			pages={709--772},
			issn={1435-9855},
			review={\MR{4859577}},
			doi={10.4171/jems/1386},
		}
		
		
		
		
		
			\bib{Patsalos-Saroglou-arxiv}{article}{
			title={A note on the $L_p$-Brunn--Minkowski inequality for intrinsic volumes and the $L_p$-Christoffel-Minkowski problem},
			author={K.  Patsalos},
			author={C. Saroglou},
			year={2025},
			eprint={2411.17896, accepted by \textit{Trans. Amer. Math. Soc.}},
			archivePrefix={arXiv},
			primaryClass={math.MG},
			url={https://arxiv.org/abs/2411.17896 (accepted by Trans. Amer. Math. Soc.)},
		}
		
		\bib{Putterman-2021}{article}{
			author={Putterman, E.},
			title={Equivalence of the local and global versions of the
				$L^p$-Brunn--Minkowski inequality},
			journal={J. Funct. Anal.},
			volume={280},
			date={2021},
			number={9},
			pages={Paper No. 108956, 20},
			issn={0022-1236},
			review={\MR{4220744}},
			doi={10.1016/j.jfa.2021.108956},
		}
		
		
		\bib{Sadovsky-ZhangGaoYong-2025}{article}{
			author={Sadovsky, S.},
			author={Zhang, G.},
			title={Brunn--Minkowski and reverse isoperimetric inequalities for dual
				quermassintegrals},
			journal={Adv. Math.},
			volume={480},
			date={2025},
			number={part A},
			part={part A},
			pages={Paper No. 110456, 14},
			issn={0001-8708},
			review={\MR{4941129}},
			doi={10.1016/j.aim.2025.110456},
		}
		
		
		
		
		\bib{book-convex-body}{book}{
			author={Schneider, R.},
			title={Convex bodies: the Brunn--Minkowski theory},
			series={Encyclopedia of Mathematics and its Applications},
			volume={151},
			edition={Second expanded edition},
			publisher={Cambridge University Press, Cambridge},
			date={2014},
			pages={xxii+736},
			isbn={978-1-107-60101-7},
			review={\MR{3155183}},
		}
		

	\bib{Sheng-Wang}{article}{
		title={On the new weighted geometric inequalities near the sphere in space forms},
		author={W. Sheng},
		author={Y. Wang},
		year={2025},
		eprint={2508.04067},
		archivePrefix={arXiv},
		primaryClass={math.DG},
		url={https://arxiv.org/abs/2508.04067},
	}
	
	
	\bib{Xi-Leng-2016}{article}{
		author={Xi, D.},
		author={Leng, G.},
		title={Dar's conjecture and the log-Brunn--Minkowski inequality},
		journal={J. Differential Geom.},
		volume={103},
		date={2016},
		number={1},
		pages={145--189},
		issn={0022-040X},
		review={\MR{3488132}},
	}
	
	
	\bib{Xi-Zhang-2022}{article}{
		author={Xi, D.},
		author={Zhang, Z.},
		title={The $L_p$ Brunn--Minkowski inequalities for dual quermassintegrals},
		journal={Proc. Amer. Math. Soc.},
		volume={150},
		date={2022},
		number={7},
		pages={3075--3086},
		issn={0002-9939},
		review={\MR{4428890}},
		doi={10.1090/proc/15952},
	}
	
	
	
	\bib{Xia-phd}{book}{
		AUTHOR = {Xia, C.},
		TITLE = {On a Class of Anisotropic Problems},
		SERIES = {PhD Thesis},
		PUBLISHER = {Albert-Ludwigs University Freiburg},
		YEAR = {2012},
	}
	
	
	
	
	
	
	
	
	\bib{Xia2017}{article}{
		title={Inverse anisotropic mean curvature flow and a Minkowski type inequality},
		author={Xia, C.},
		JOURNAL = {Adv. Math.},
		FJOURNAL = {Advances in Mathematics},
		volume={315},
		pages={102-129},
		year={2017},
	}
	
	
	
	
\end{biblist}
\end{document}